\documentclass[12pt]{amsart}
\usepackage{amsmath,amssymb,amsthm}
\usepackage[margin=1in]{geometry}
\usepackage{cleveref}
\usepackage{setspace}
\newtheorem{thm}{Theorem}[section]
\newtheorem{prop}[thm]{Proposition}
\newtheorem{cor}[thm]{Corollary}
\newtheorem{lem}[thm]{Lemma}
\newtheorem{defn}[thm]{Definition}
\newtheorem{exam}[thm]{Example}
\theoremstyle{remark}
\newtheorem{rmk}[thm]{Remark}

\begin{document}

\title{A Scalar-Valued Weyl Transform Associated to the Heisenberg Group}
\author{R. Lakshmi Lavanya}
\address{Department of Mathematics, Indian Institute of Science Education and Research (IISER) Tirupati, Yerpedu, Tirupati 517619, Andhra Pradesh, India.}
\email{rlakshmilavanya@labs.iisertirupati.ac.in}
\date{}

\begin{abstract}
The classical group Fourier transform on the Heisenberg group factors through the ordinary Weyl transform. Motivated by this, we ask
what plays the same role for the Strichartz Fourier transform, a scalar-valued Fourier transform on the Heisenberg group,
recently introduced by Thangavelu. We identify an operator, that factorises on each ray of the Heisenberg
fan, and we call it the Strichartz Weyl transform. Up to a scalar,
this turns out to be a spectral projection of the special
Hermite operator, and so inherits a mature body of
independently developed spectral theory. We develop the fixed-ray
structural theory of the Strichartz Weyl transform, including boundedness,
idempotency, self-adjointness, a self-reproducing range condition,
geometric covariance, and a Plancherel identity, and use it for two
purposes. For the Strichartz Weyl transform, we obtain a Fock-space transport via
the Bargmann transform, and a single-ray Benedicks-type uncertainty
theorem from the ellipticity of the special Hermite operator. Then we assemble this theory across the Heisenberg fan to recover
the standard foundational toolkit any new transform requires, such as 
inversion, Plancherel, a unitary-range characterisation, a
Hausdorff-Young inequality, and Riemann-Lebesgue and continuity theory
for the Strichartz Fourier transform, largely for free, and give a Bochner-Riesz summability theorem from
inversion. We also
treat the Hecke-Bochner formula, and discuss the image of
Schwartz space, under the Strichartz Weyl transform, established at a single fixed ray. Finally, we give a detailed account of two results beyond
the reach of the fixed-ray method, Ingham's theorem and a joint-ray
Benedicks theorem recently obtained by Dabra, Dasgupta, and Gulia, and
explain precisely why.
\end{abstract}

\subjclass[2020]{Primary 43A85, 42C05; Secondary 33C45, 35P10}

\keywords{Heisenberg group, Strichartz Fourier transform, Weyl transform, special Hermite expansions, spectral theory, Gelfand transform, Hecke-Bochner formula, Bochner-Riesz means, uncertainty principle}

\maketitle

\section{Introduction}
The group Fourier transform $\widehat f(\lambda)$ on the Heisenberg
group $\mathbb{H}^n$ is well studied
\cite{Folland1989,Thangavelu1998}. It is operator-valued, taking
values in $B(L^2(\mathbb{R}^n))$, which makes it more challenging to
work with for some purposes, such as describing the image of the
Schwartz space or formulating sharp weighted inequalities. Thangavelu
recently introduced a scalar-valued alternative on $\mathbb{H}^n$, the
Strichartz Fourier transform $\mathcal{F}_S$, built from the joint
spectral theory of the sublaplacian $\mathcal{L}$ and the central
vector field $T = \partial/\partial t$ \cite{Thangavelu2023}.

For the classical group Fourier transform, the factorisation
\begin{equation} \label{eq:classical-factorization}
\widehat f(\lambda) = W_\lambda(f^\lambda)
\end{equation}
is more than just a structural identity, where $W_\lambda$ denotes the
classical Weyl transform (recalled in \Cref{sec:prelims}). So far, it
has been the route to many results on the group Fourier transform. A
fact about $\widehat f(\lambda)$ is typically established by first
proving the corresponding fact for $W_\lambda$, an operator on
$\mathbb{C}^n$ that is usually more tractable than $\widehat
f(\lambda)$ itself. The result is then transported across the
factorisation. This method works precisely because $\widehat f(\lambda)$ and
$W_\lambda$, despite being defined independently of one another, are
tied together by \eqref{eq:classical-factorization}.

The definition of $\mathcal{F}_S$ in \cite{Thangavelu2023} already
carries the raw ingredients of an analogous relation. On each ray of
$\Omega$, $\mathcal{F}_S(f)(a,z)$ is built from $f^{-\lambda}$ (the
partial Fourier transform of $f$ in the central variable), composed with twisted
convolution against the Laguerre kernel $\varphi^{n-1}_{k,\lambda}$.
This is the same shape as the classical case, where $\widehat
f(\lambda)$ is built from $f^\lambda$ composed with integration
against $\pi_\lambda(z,0)$. This raises a natural question - does $\mathcal{F}_S$ admit a
factorisation of the same kind, through some companion operator on
$\mathbb{C}^n$? If so, the same fork, establishing the fact for the
companion operator and then transporting it, would open up for
theorems on $\mathcal{F}_S$.

This paper answers this affirmatively. We name this
new transform the Strichartz Weyl transform and denote it by
$W_{S,\lambda,k}$ (cf.~\Cref{def:swtdef}). When $k$ is understood, we
write $W_{S,\lambda}$ in its place. This gives the following factorisation for
$\mathcal{F}_S$:
\begin{equation}
\label{eq:sftfactor}
\mathcal{F}_S(f)(a,z) = W_{S,\lambda}\big(f^{-\lambda}\big)(z), \qquad
a = (\lambda, (2k+n)|\lambda|) \in R_k.
\end{equation}
A priori, the operator completing this analogy need not have been
familiar. However, up to scaling, it proves to be one, as
it is the $k$-th spectral projection of the special Hermite operator. This is the natural analogue of the Laplacian governing twisted convolution
on $\mathbb{C}^n$. This identification is a second and independently
notable fact (cf.~\Cref{prop:swt-projection}). This means
$\mathcal{F}_S$ inherits, through the factorisation, the spectral
theory already developed for this projection elsewhere, independently
of $\mathcal{F}_S$, in the work of Peetre \cite{Peetre1972}, Koch and
Ricci \cite{KochRicci2007}, and Thangavelu
\cite{ThangaveluLectures1993}.

This is the same pattern by which \cite{Thangavelu2023} itself arose: the joint eigenfunctions $e^{n-1}_{k,\lambda}$ and the Laguerre-kernel twisted convolution were already present in the literature on special Hermite expansions, but had not been isolated as the object defining a scalar-valued transform in its own right until \cite{Thangavelu2023} did so. The operator $W_{S,\lambda,k}$ we isolate here is, in the same sense, already written into ~\Cref{def:swtdef}. The content of \Cref{prop:swt-projection} is not that this fact is deep. Rather, naming the operator and developing it independently of $\mathcal{F}_S$ makes available, for free, a body of structural theory that stays hidden if one looks only at~\Cref{def:swtdef}.

We note that not every structural fact we establish about
$W_{S,\lambda}$ is intended, or needs, to yield a correspondingly
simplified proof of some fact about $\mathcal{F}_S$. The classical Weyl
transform $W_\lambda$ has long been studied as an object in its own
right, and results about it are not routinely justified by whether
they shorten some downstream statement about $\widehat f(\lambda)$. We
extend the same treatment to $W_{S,\lambda}$. Some of the structural
theory developed in \Cref{sec:swt} and some of the fixed-ray results
of \Cref{sec:core-toolkit} are used below to recover theorems for
$\mathcal{F}_S$. Some results are recorded as properties of $W_{S,\lambda}$
in their own right, in keeping with how the classical theory has
always treated $W_\lambda$.

\paragraph{\textbf{Organisation of the paper.}}
The aim of the present paper is to put this inheritance to systematic
use, and the value of doing so lies in its reach, not in the depth of any single derivation. The point of recognising
$W_{S,\lambda}$ as a scaled spectral projection is that a whole body of
structural facts about such projections, developed independently of
$\mathcal{F}_S$ and in some cases decades before it, becomes available
essentially for free. The paper is organised around this distinction.\\
\textbf{\Cref{sec:prelims}} recalls definitions and fixes notations to be followed in the subsequent sections.\\
\textbf{\Cref{sec:swt}} develops $W_{S,\lambda}$ at a single fixed ray as a self-contained object, and proves boundedness, symmetries, idempotency,
self-adjointness, a self-reproducing range condition, geometric
covariance, and the Plancherel bound.\\
\textbf{\Cref{sec:new-results}} is where the dictionary is put to work to obtain new results for the Strichartz Fourier transform. We obtain a Fock-space transport of $W_{S,\lambda}$ via
the Bargmann transform, and a Benedicks-type uncertainty theorem obtained
directly from the ellipticity of the special Hermite operator. These demonstrate the value of the Strichartz Weyl transform as a tool for obtaining new results about the Strichartz Fourier transform.\\
\textbf{\Cref{sec:core-toolkit}} assembles this fixed-ray theory across $\Omega$ to recover the
standard foundational toolkit that any new transform requires, namely,
inversion, Plancherel, a unitary-range characterisation, a
Hausdorff-Young inequality, and the Riemann-Lebesgue and continuity
theory. Each of these was already established for $\mathcal{F}_S$ in
\cite{Thangavelu2023} by direct computation with Laguerre and Bessel
asymptotics. We include this material not because it is new, but
because its near-triviality once routed through $W_{S,\lambda}$ is
itself the clearest demonstration of the dictionary's value.\\
\textbf{\Cref{sec:hecke-bochner-schwartz}} treats the Hecke-Bochner formula and the image of $\mathcal{S}(\mathbb{H}^n)$ under $\mathcal{F}_S$. This material is of a different character from \Cref{sec:core-toolkit,sec:new-results}. The Hecke-Bochner formula itself is reached through $W_{S,\lambda}$, but the full Schwartz-space image is not. It requires a compatibility condition linking values on different rays as $\lambda \to 0$, which is invisible to $W_{S,\lambda}$ at any single fixed ray, and we point instead to the harder global direction of a theorem of Astengo, Di Blasio, and Ricci \cite{AstengoDiBlasioRicci2007}, as restated by Thangavelu \cite{Thangavelu2023}, for that part.\\
\textbf{\Cref{sec:limits}} gives a full account of two results
that lie beyond this method altogether. We discuss the Ingham's uncertainty theorem,
and the genuinely joint-ray Benedicks-type theorem recently obtained
by Dabra, Dasgupta, and Gulia \cite{DabraDasguptaGulia2026}. Both are negative results about the fixed-ray method, each for a specific reason, and we discuss both as carefully as the positive results, rather than mentioning them only in passing.

\section{Preliminaries}
\label{sec:prelims}
In this section, we recall the definitions and notation used throughout the paper, covering the Heisenberg group and its representations, the Laguerre functions and special Hermite operator, the Heisenberg motion group, and the Heisenberg fan on which the Strichartz Fourier transform is defined.

\subsection{The Heisenberg group, representations, and twisted convolution}
Let $\mathbb{H}^n = \mathbb{C}^n \times \mathbb{R}$ be the
$(2n+1)$-dimensional Heisenberg group with group law
\begin{equation*}
(z,t)(w,s) = \Big(z+w,\ t+s+\tfrac{1}{2}\operatorname{Im}(z\cdot\bar w)\Big).
\end{equation*}
The Lebesgue measure $dz\,dt$ serves as Haar measure on $\mathbb{H}^n$. For $\lambda \in
\mathbb{R}^*$, write
\begin{equation*}
f^\lambda(z) := \int_{\mathbb{R}} f(z,t)\, e^{i\lambda t}\, dt, \ z\in \mathbb{C}^n.
\end{equation*}

For $\lambda \in \mathbb{R}^*$, the Schr\"odinger representation
$\pi_\lambda$ acts on $L^2(\mathbb{R}^n)$ by
\begin{equation*}
\pi_\lambda(z,t)\phi(\xi) = e^{i\lambda t}\, e^{i\lambda(x\cdot\xi +
\frac{1}{2}x\cdot y)}\, \phi(\xi+y), \qquad z = x+iy.
\end{equation*}
These representations are unitary and irreducible, and by the
Stone-von Neumann theorem account, up to unitary equivalence, for all
infinite-dimensional irreducible unitary representations of
$\mathbb{H}^n$ acting as $e^{i\lambda t}I$ on the centre.

For $f \in L^1(\mathbb{H}^n)$, the group Fourier transform of $f$ is the operator on
$L^2(\mathbb{R}^n)$ given by
\begin{equation*}
\widehat f(\lambda) := \int_{\mathbb{H}^n} f(z,t)\, \pi_\lambda(z,t)\, dz\,dt, \ \lambda \in \mathbb{R}^*.
\end{equation*}
This operator is bounded, with $\|\widehat f(\lambda)\|_{\mathrm{op}} \le \|f\|_1$.

The classical Weyl transform $W_\lambda$, $\lambda \in \mathbb{R}^*$, acts on
$L^2(\mathbb{R}^n)$ by
\begin{equation*}
W_\lambda(g) := \int_{\mathbb{C}^n} g(z)\, \pi_\lambda(z,0)\, dz, \qquad g \in L^1(\mathbb{C}^n).
\end{equation*}
$W_\lambda$ has long been studied as an object in its own right, independently of the
group Fourier transform.

The group Fourier transform of $f \in L^1(\mathbb{H}^n)$ factors through the Weyl
transform as
\begin{equation}
\widehat f(\lambda) = W_\lambda(f^\lambda).
\end{equation}

For $\lambda \in \mathbb{R}^*$, the $\lambda$-twisted convolution of
$f,g \in L^1(\mathbb{C}^n)$ is
\begin{equation*}
(f *_\lambda g)(z) = \int_{\mathbb{C}^n} f(z-w)\, g(w)\,
e^{i\frac{\lambda}{2}\operatorname{Im}(z\cdot\bar w)}\, dw.
\end{equation*}
This product is associative on $L^1(\mathbb{C}^n)$, and the
substitution $w \mapsto z-w$ gives the identity $f *_\mu g = g *_{-\mu}
f$, $\mu \in \mathbb{R}^*$. Since complex conjugation of $\mathbb{C}^n$ preserves radial
functions and reverses the sign of $\operatorname{Im}(z\cdot\bar w)$,
twisted convolution of radial functions is commutative and does not
depend on the sign of the parameter \cite{ThangaveluLectures1993}.

\subsection{Laguerre functions and the special Hermite operator}

\label{sec:laguerre-hermite}

For $\delta \ge -1/2$, let $L^\delta_k$ denote the Laguerre polynomials
of type $\delta$. For $\lambda \in \mathbb{R}^*$ and $k \in
\mathbb{N}$, the Laguerre functions are
\begin{equation*}
\varphi^{n-1}_{k,\lambda}(z) = L^{n-1}_k\!\left(\tfrac{1}{2}|\lambda||z|^2\right)
e^{-\frac{1}{4}|\lambda||z|^2}.
\end{equation*}
These satisfy the uniform bound
\begin{equation*}
|\varphi^{n-1}_{k,\lambda}(z)| \le L^{n-1}_k(0), \qquad z \in \mathbb{C}^n,
\end{equation*}
cf.\ \cite[\S 7.21]{Szego1975}. The normalised Laguerre kernel is
\begin{equation*}
\widetilde\varphi^{n-1}_{k,\lambda} := \frac{\varphi^{n-1}_{k,\lambda}}{L^{n-1}_k(0)}
= \frac{k!\,(n-1)!}{(k+n-1)!}\,\varphi^{n-1}_{k,\lambda},
\end{equation*}
satisfying $|\widetilde\varphi^{n-1}_{k,\lambda}| \le 1$ and
$\widetilde\varphi^{n-1}_{k,\lambda}(0) = 1$. Writing
\begin{equation*}
e^{n-1}_{k,\lambda}(z,t) := e^{i\lambda t}\, \varphi^{n-1}_{k,\lambda}(z)
\end{equation*}
gives the joint eigenfunctions of the sublaplacian $\mathcal{L}$ and $T
= \partial/\partial t$ on $\mathbb{H}^n$, satisfying \[
\mathcal{L} e^{n-1}_{k,\lambda} = (2k+n)|\lambda|\, e^{n-1}_{k,\lambda}, \qquad
T e^{n-1}_{k,\lambda} = i \lambda\, e^{n-1}_{k,\lambda}.
\]

The special Hermite operator $L_{-\lambda}$ has $L^2(\mathbb{C}^n)$-eigenspaces
indexed by $j \in \mathbb{N}$, corresponding to eigenvalue
$(2j+n)|\lambda|$. Let $P_j$ denote the orthogonal projection onto the
$j$-th eigenspace. When the dependence on $\lambda$ needs to be emphasised,
we write $P_j(\lambda) := P_j$. Then
\begin{equation*}
u = \sum_{j=0}^\infty P_j u, \qquad u \in L^2(\mathbb{C}^n),
\end{equation*}
and
\begin{equation}\label{eq:laguerre-projection}
u *_{-\lambda} \varphi^{n-1}_{j,\lambda} = (2\pi)^n |\lambda|^{-n} P_j u.
\end{equation}
The radial functions in the $j$-th eigenspace are the multiples of
$\varphi^{n-1}_{j,\lambda}$. In particular, for $i,j \in \mathbb{N}$,
\begin{equation*}
\varphi^{n-1}_{i,\lambda} *_{-\lambda} \varphi^{n-1}_{j,\lambda}
= \delta_{ij}\,(2\pi)^n |\lambda|^{-n}\, \varphi^{n-1}_{j,\lambda}.
\end{equation*}

\begin{rmk}
The Laguerre kernel $\varphi^{n-1}_{k,\lambda}$ depends on $\lambda$
only through $|\lambda|$, so $\varphi^{n-1}_{k,\lambda} =
\varphi^{n-1}_{k,-\lambda}$. The special Hermite projection
formula in \cite{Thangavelu2023} is stated for $L_\lambda$, via
$+\lambda$-twisted convolution, as
\begin{equation*}
g *_\lambda \varphi^{n-1}_{k,\lambda} = (2\pi)^n|\lambda|^{-n}\, P_k(\lambda)\, g,
\qquad g \in L^2(\mathbb{C}^n).
\end{equation*}
Relabelling $\lambda$ as $-\lambda$ on both sides, and using
$\varphi^{n-1}_{k,-\lambda} = \varphi^{n-1}_{k,\lambda}$, reproduces
the identity above verbatim,
\begin{equation*}
g *_{-\lambda} \varphi^{n-1}_{k,\lambda} = (2\pi)^n|\lambda|^{-n}\, P_k(-\lambda)\, g.
\end{equation*}
Here $P_j$ is now understood as the projection onto the $j$-th eigenspace
of $L_{-\lambda}$ rather than $L_\lambda$. We adopt this $L_{-\lambda}$, that is, the $-\lambda$-convention throughout, since it matches the
defining formula for $\mathcal{F}_S$ directly
(cf. Remark~\ref{rmk:comparing-FS}).
\end{rmk}

Let $\{\Phi^\lambda_\alpha\}_{\alpha\in\mathbb{N}^n}$ denote the scaled real Hermite functions on $\mathbb{R}^n$, given by 
\[
\Phi^\lambda_\alpha(x) := |\lambda|^{n/4}\Phi_\alpha(\sqrt{|\lambda|}x).
\]
The family $\{\Phi^\lambda_\alpha\}_{\alpha\in\mathbb{N}^n}$ forms an orthonormal basis of $L^2(\mathbb{R}^n)$. For $\alpha, \beta \in \mathbb{N}^n$, the \emph{special Hermite functions} $\Phi^\lambda_{\alpha\beta}$ are defined by
\[
\Phi^\lambda_{\alpha\beta}(z) := (2\pi)^{-n/2}\big(\pi_{-\lambda}(z,0)\Phi^\lambda_\alpha,\ \Phi^\lambda_\beta\big).
\]
These satisfy
\begin{equation*}
\Phi^\lambda_{\alpha\beta} *_{-\lambda} \Phi^\lambda_{\mu\nu}
= (2\pi)^{n/2}|\lambda|^{-n}\, \delta_{\beta\mu}\, \Phi^\lambda_{\alpha\nu},
\end{equation*}
and
\begin{equation*}
\varphi^{n-1}_{j,\lambda} = (2\pi)^{n/2} \sum_{|\gamma|=j} \Phi^\lambda_{\gamma\gamma}.
\end{equation*}
The first identity follows from the Schur
(Moyal) orthogonality relations for matrix coefficients of
$\pi_{-\lambda}$ \cite{ThangaveluLectures1993}, together with
$\pi_{-\lambda}(w,0)\pi_{-\lambda}(z-w,0) =
\pi_{-\lambda}(z,0)e^{i\lambda/2\,\operatorname{Im}(w\bar z)}$.
The next one records that the sum of the diagonal matrix
coefficients over the $j$-th eigenspace of the harmonic oscillator is
the Laguerre function.  Together,
these two identities give
\begin{equation*}
\operatorname{Ran}(P_k) = \overline{\operatorname{span}}
\{\Phi^\lambda_{\alpha\beta} : |\beta| = k\}.
\end{equation*}

\subsection{The Heisenberg motion group and class-one representations}

\label{sec:motiongroup}

The unitary group $U(n)$ acts on $\mathbb{H}^n$ by automorphisms,
\[
\sigma.(z,t) = (\sigma z,t),
\]
and $G_n := \mathbb{H}^n \rtimes U(n)$ denotes the Heisenberg motion
group, with group law
\begin{equation*}
(z,t,\sigma)(w,s,\tau) = \Big(z+\sigma w,\ t+s+\tfrac{1}{2}
\operatorname{Im}(z\cdot\overline{\sigma w}),\ \sigma\tau\Big).
\end{equation*}
The Haar measure on $G_n$ is $d\sigma\,dz\,dt$, with $d\sigma$ the
normalised Haar measure on $U(n)$. Functions on $\mathbb{H}^n$ are
precisely the right $U(n)$-invariant functions on $G_n$. For
$\alpha,\beta\in\mathbb{N}^n$, let
\begin{equation*}
E^\lambda_{\alpha,\beta}(z,t) := \big(\pi_\lambda(z,t)\Phi^\lambda_\alpha,\ \Phi^\lambda_\beta\big)
\end{equation*}
denote the matrix coefficients of the Schr\"odinger representation.
For $\lambda \ne 0$ and $k \in \mathbb{N}$, let $\mathcal{H}^\lambda_k$
denote the Hilbert space spanned by $\{E^\lambda_{\alpha,\beta} :
\alpha, \beta \in \mathbb{N}^n,\ |\beta|=k\}$, equipped with the inner
product
\begin{equation*}
(f,g)_{\mathcal{H}^\lambda_k} := (2\pi)^{-n}|\lambda|^n
\int_{\mathbb{C}^n} f(z,0)\, \overline{g(z,0)}\, dz.
\end{equation*}
For $\lambda \ne 0$ and $k \in \mathbb{N}$, let $A_k^{(\lambda)}$
denote the orthogonal projection of $L^2(\mathbb{R}^n)$ onto the span
of $\{\Phi^\lambda_\alpha : \alpha \in \mathbb{N}^n,\ |\alpha|=k\}$,
the $k$-th eigenprojection of the scaled Hermite operator. This corresponds,
under the usual identification, to the projection denoted $P_k$ in
\Cref{sec:laguerre-hermite} above, and $P_k(\lambda)$ in \cite{Thangavelu2023}. The representation
$\rho^\lambda_k$ of $G_n$ on $\mathcal{H}^\lambda_k$ is defined by
\begin{equation*}
\rho^\lambda_k(z,t,\sigma)\phi(w,s) := \phi\big((z,t,\sigma)^{-1}(w,s)\big),
\qquad (w,s) \in \mathbb{H}^n.
\end{equation*}
This is an irreducible unitary representation of $G_n$ for all
$\lambda \ne 0$, $k \in \mathbb{N}$, and $e^{n-1}_{k,\lambda}$ is, up
to scalar, its unique $U(n)$-fixed vector. So $\rho^\lambda_k$ is a
class-one representation of the pair $(G_n,U(n))$. This indexing of
eigenspace membership of $E^\lambda_{\alpha,\beta}$ and
$\Phi^\lambda_{\alpha\beta}$, determined by the second index, is the
convention adopted throughout.

\begin{prop}[Motion-group realisation, see e.g.~\cite{Thangavelu1998}]
For $f \in L^1(\mathbb{H}^n)$, viewed as a right $U(n)$-invariant
function on $G_n$,
\begin{equation*}
\widehat f(a,z) = \rho^\lambda_k(f)\, e_a(z,0), \qquad \rho^\lambda_k(f)
:= \int_{G_n} f(g)\, \rho^\lambda_k(g)\, dg,
\end{equation*}
where $e_a := e^{n-1}_{k,\lambda}$ for $a = (\lambda,(2k+n)|\lambda|)$.
\end{prop}
This is the motion-group realisation of the Strichartz Fourier transform.

\begin{rmk}
The classification of bounded spherical functions for $(G_n,U(n))$
underlying the identification of $\Omega$ with its Gelfand spectrum
(\Cref{sec:gelfand} below) is due to \cite{BensonJenkinsRatcliff1998}. A further classification of bounded $K$-spherical functions on
Heisenberg groups appears in \cite{BensonJenkinsRatcliff1992}. Each
projection $A_k^{(\lambda)}$ is the orthogonal projection onto the
$k$-th $U(n)$-isotypic subspace of $\pi_\lambda|_{U(n)}$. This is the
representation-theoretic reason that radial, that is,
$U(n)$-invariant, functions on $\mathbb{H}^n$ correspond to scalar
spectral values indexed by $\Omega$.
\end{rmk}

\subsection{The Heisenberg fan and the Strichartz Fourier transform}

The Heisenberg fan is
\begin{equation*}
\Omega = \bigcup_{k=0}^\infty R_k \cup R_\infty, \qquad R_k =
\{(\lambda,(2k+n)|\lambda|) : \lambda \in \mathbb{R}^*\}, \qquad R_\infty
= \{(0,\tau) : \tau \ge 0\},
\end{equation*}
topologised so that $a_j = (\lambda_j,k_j) \to (0,\tau)$ whenever
$\lambda_j \to 0$ and $(2k_j+n)|\lambda_j| \to \tau$.\\
For $a =
(\lambda,(2k+n)|\lambda|) \in R_k$ write $e_a := e^{n-1}_{k,\lambda}$.
As a subset of $\mathbb{R}^2$, $\Omega$ inherits the Euclidean metric, and with it the subspace topology. On $\Omega$, the measure $\nu$ is defined by
\begin{equation*}
\int_\Omega \phi(a)\, d\nu(a) := (2\pi)^{-2n-1} \int_{-\infty}^\infty
\sum_{k=0}^\infty \phi\big(\lambda,(2k+n)|\lambda|\big)\, |\lambda|^{2n}\,
d\lambda, \qquad \phi \in C_c(\Omega).
\end{equation*}

\begin{defn}(\cite{Thangavelu2023})
\label{def:sftdef}
For $f \in L^1 \cap L^2(\mathbb{H}^n)$ and $a = (\lambda,(2k+n)|\lambda|) \in R_k$,
the Strichartz Fourier transform is defined by
\begin{equation*}
\mathcal{F}_S(f)(a,z) =\widehat{f}(a,z) := f^{-\lambda} *_{-\lambda} \varphi^{n-1}_{k,\lambda}(z),
\end{equation*}
and, for $a = (0,\tau) \in R_\infty$,
\begin{equation*}
\mathcal{F}_S(f)(0,\tau,z) = \widehat{f}(0,\tau,z) := (n-1)!\, 2^{n-1} \int_{\mathbb{H}^n} f(w,s)\,
\frac{J_{n-1}(\sqrt{\tau}|z-w|)}{(\sqrt{\tau}|z-w|)^{n-1}}\, dw\, ds,
\end{equation*}
where $J_{n-1}$ is the Bessel function of order $n-1$, and
\[
\mathcal{F}_S(f)(0,0,z) = \widehat f(0,0,z) := \int_{\mathbb{H}^n} f(w,s)\, dw\, ds.
\]
\end{defn}

The normalised Strichartz Fourier transform
is
\begin{align*}
\widetilde f(a,z) &:= \frac{k!(n-1)!}{(k+n-1)!}\, \widehat f(a,z), \quad (a,z)\in R_k, \\
\widetilde f(0,\tau,z) &:= \widehat f(0,\tau,z), \quad (0,\tau) \in R_\infty.
\end{align*}

Using the identity $\varphi^{n-1}_{k,\lambda} *_\lambda
\varphi^{n-1}_{k,\lambda}(z) = (2\pi)^n|\lambda|^{-n}\varphi^{n-1}_{k,\lambda}(z)$,
the Fourier transform $\widehat f(a,z)$ satisfies, for $a \in R_k$,
\begin{equation*}
(2\pi)^{-n}|\lambda|^n \int_{\mathbb{C}^n} \widehat f(a,w)\,
\varphi^{n-1}_{k,\lambda}(z-w)\, e^{i\frac{\lambda}{2}\operatorname{Im}(z\cdot\bar w)}\, dw
= \widehat f(a,z).
\end{equation*}
If $L^2_0(\Omega\times\mathbb{C}^n, d\nu\, dw)$ denotes the subspace of
$L^2(\Omega\times\mathbb{C}^n, d\nu\, dw)$ consisting of functions
$F(a,z)$ satisfying the above condition at every $a \in R_k$, then
$\mathcal{F}_S(f) \in L^2_0(\Omega\times\mathbb{C}^n, d\nu\,dw)$
whenever $f \in L^1 \cap L^2(\mathbb{H}^n)$.

\section{The Strichartz Weyl Transform}
\label{sec:swt}

This section defines the Strichartz Weyl transform and studies some of its basic properties at a single fixed ray.

\begin{defn}
\label{def:swtdef}
For $\lambda \in \mathbb{R}^*$ and $k \in \mathbb{N}$, the Strichartz Weyl
transform is the operator
\begin{equation*}
W_{S,\lambda,k}(g) := g *_{-\lambda} \varphi^{n-1}_{k,\lambda}, \qquad
g \in L^1 \cap L^2(\mathbb{C}^n),
\end{equation*}
matching the domain on which $\mathcal{F}_S$ itself is defined
(\Cref{def:sftdef}). When the index $k$ is understood, we write
$W_{S,\lambda}(g)$ in place of $W_{S,\lambda,k}(g)$.
\end{defn}
Let \[
\widetilde{W}_{S,\lambda,k}(g) := g *_{-\lambda} \widetilde{\varphi}^{\,n-1}_{k,\lambda},\qquad
g \in L^1 \cap L^2(\mathbb{C}^n),
\] denote the \textit{normalised
Strichartz Weyl transform.} 

\begin{prop}
\label{prop:swt-projection}
For $\lambda \in \mathbb{R}^*$ and $k \in \mathbb{N}$, let $P_k$ denote
the $k$-th spectral projection of $L_{-\lambda}$ on $L^2(\mathbb{C}^n)$.
Then, for every $g \in L^1 \cap L^2(\mathbb{C}^n)$,
\begin{equation*}
W_{S,\lambda,k}(g) = (2\pi)^n |\lambda|^{-n} P_k g.
\end{equation*}
\end{prop}

\begin{proof}
Recall from \Cref{sec:laguerre-hermite} that the special Hermite operator $L_{-\lambda}$
has $L^2(\mathbb{C}^n)$-eigenspaces indexed by $j \in \mathbb{N}$,
corresponding to eigenvalue $(2j+n)|\lambda|$, with $P_j$ the orthogonal
projection onto the $j$-th eigenspace. These projections are
realised through twisted convolution against the Laguerre kernel via \eqref{eq:laguerre-projection}. Taking $j = k$ and $u = g \in L^1 \cap L^2(\mathbb{C}^n)$ gives
\begin{equation*}
g *_{-\lambda} \varphi^{n-1}_{k,\lambda} = (2\pi)^n |\lambda|^{-n} P_k g,
\end{equation*}
and the left-hand side is $W_{S,\lambda,k}(g)$ by \Cref{def:swtdef}. As
both sides are defined for every $g \in L^1 \cap L^2(\mathbb{C}^n)$,
this proves the stated identity on this common domain.
\end{proof}

\begin{rmk}
\label{rmk:comparing-FS}
This is the factorisation of $\mathcal{F}_S$ promised in the introduction. Comparing
\Cref{def:swtdef} with the definition of $\mathcal{F}_S$ (\Cref{def:sftdef}) gives
\begin{equation*}
\mathcal{F}_S(f)(a,z) = W_{S,\lambda}\big(f^{-\lambda}\big)(z), \qquad
a = (\lambda,(2k+n)|\lambda|) \in R_k,
\end{equation*}
for $f \in L^1 \cap L^2(\mathbb{H}^n)$, since $f^{-\lambda} \in
L^1\cap L^2(\mathbb{C}^n)$ for almost every $\lambda$ whenever
$f \in L^1\cap L^2(\mathbb{H}^n)$.
\end{rmk}

\subsection{Tautological properties}
The following properties follow immediately from \Cref{prop:swt-projection}
and the corresponding properties of $P_k$, with no further argument required.
Since $P_k^2 = P_k$,
\begin{equation*}
W_{S,\lambda} \circ W_{S,\lambda} = (2\pi)^n |\lambda|^{-n}\, W_{S,\lambda}.
\end{equation*}
This means $W_{S,\lambda}$ is idempotent up to the scalar $(2\pi)^n|\lambda|^{-n}$,
the same scalar appearing in \Cref{prop:swt-projection}.\\
Since $P_k^* = P_k$,
\begin{equation*}
(W_{S,\lambda}u,v) = (u,W_{S,\lambda}v), \qquad u,v \in L^2(\mathbb{C}^n),
\end{equation*}
that is, $W_{S,\lambda}$ is self-adjoint on $L^2(\mathbb{C}^n)$.
\subsection{Boundedness and symmetries}

\begin{prop}
For $g \in L^1(\mathbb{C}^n)$,
\begin{equation*}
\|W_{S,\lambda}(g)\|_\infty \le L^{n-1}_k(0)\, \|g\|_1.
\end{equation*}
\end{prop}

\begin{proof}
Since $|\varphi^{n-1}_{k,\lambda}| \le L^{n-1}_k(0)$, we get
\begin{equation*}
|W_{S,\lambda}(g)(z)| \le \|g\|_1 \  \sup_{\mathbb{C}^n} | \ \varphi^{n-1}_{k,\lambda}|
= \|g\|_1\, L^{n-1}_k(0).
\end{equation*}
\end{proof}

\begin{prop}[Reflection]
For $g \in L^1(\mathbb{C}^n)$ and $z \in \mathbb{C}^n$, writing
$g^*(z) := g(-z)$,
\begin{equation*}
W_{S,\lambda}(g)(-z) = W_{S,\lambda}(g^*)(z).
\end{equation*}
\end{prop}

\begin{proof}
Substituting $y \mapsto -y$ in the definition of $W_{S,\lambda}(g)(-z)$
and using evenness of $\varphi^{n-1}_{k,\lambda}$ gives the stated
identity.
\end{proof}

\begin{prop}[Conjugation]
\label{prop:conjugation}
For $g \in L^1(\mathbb{C}^n)$,
\begin{equation*}
W_{S,-\lambda}(\bar g) = \overline{W_{S,\lambda}(g)}.
\end{equation*}
\end{prop}

\begin{proof}
By definition,
\[
W_{S,-\lambda}(\bar g)(z) = \int_{\mathbb{C}^n} \bar g(z-w)\,\varphi^{n-1}_{k,-\lambda}(w)\,
e^{i\frac{\lambda}{2}\mathrm{Im}(z\cdot\bar w)}\, dw.
\]
Since $\varphi^{n-1}_{k,\lambda}$ depends on $\lambda$ only through
$|\lambda|$ and is real-valued (\Cref{sec:laguerre-hermite}), we have \\
$\varphi^{n-1}_{k,-\lambda} = \varphi^{n-1}_{k,\lambda}$. This gives
\begin{align*}
W_{S,-\lambda}(\bar g)(z) &= \int_{\mathbb{C}^n} \overline{g(z-w)}\,\varphi^{n-1}_{k,\lambda}(w)\,
e^{i\frac{\lambda}{2}\mathrm{Im}(z\cdot\bar w)}\, dw \\
&= \overline{W_{S,\lambda}(g)(z)},
\end{align*}
using that $\varphi^{n-1}_{k,\lambda}$ and the measure $dw$ are real.
This proves the identity.
\end{proof}

\subsection{Radial convolution, kernel normalisation, and geometric covariance}

\begin{prop}
$W_{S,\lambda}$ satisfies the following.
\begin{enumerate}
\item[(i)] For every $u \in L^1(\mathbb{C}^n)$ and every $U(n)$-radial
$v \in L^1(\mathbb{C}^n)$,
\begin{equation*}
W_{S,\lambda}(u *_{-\lambda} v) = \frac{W_{S,\lambda}(v)(0)}{L^{n-1}_k(0)}\, W_{S,\lambda}(u).
\end{equation*}
\item[(ii)] $W_{S,\lambda}\big(\varphi^{n-1}_{k,\lambda}\big) =
(2\pi)^n |\lambda|^{-n}\, \varphi^{n-1}_{k,\lambda}$.
\item[(iii)] For every $c \in \mathbb{C}^n$ and $g \in L^1(\mathbb{C}^n)$,
\begin{equation*}
W_{S,\lambda}(\tau_c g)(w) = e^{i\frac{\lambda}{2}\operatorname{Im}(c\cdot\bar w)}\,
W_{S,\lambda}(M_c g)(w+c),
\end{equation*}
where $\tau_c g(x) := e^{i\frac{\lambda}{2}\operatorname{Im}(x\cdot\bar c)} g(x+c)$
and $M_c g(x) := e^{i\lambda\operatorname{Im}(x\cdot\bar c)} g(x)$.
\item[(iv)] For every $\sigma \in U(n)$ and $g \in L^1(\mathbb{C}^n)$,
\begin{equation*}
W_{S,\lambda}(\rho_\sigma g)(z) = W_{S,\lambda}(g)(\sigma z),
\end{equation*}
where $\rho_\sigma g(x) := g(\sigma x)$.
\end{enumerate}
\end{prop}

\begin{proof}
\textit{(i)} The function $v$ is radial, and $P_k$ commutes with the $U(n)$-action. Hence
\begin{align*}
v *_{-\lambda} \varphi^{n-1}_{k,\lambda} = (2\pi)^n|\lambda|^{-n}P_kv
\end{align*}
is radial and lies in the range of $P_k$. The radial part of this eigenspace is one-dimensional, spanned by $\varphi^{n-1}_{k,\lambda}$. So $v *_{-\lambda} \varphi^{n-1}_{k,\lambda}$ is a scalar multiple of $\varphi^{n-1}_{k,\lambda}$. Evaluating at $z=0$ identifies the scalar as $W_{S,\lambda}(v)(0)/L^{n-1}_k(0)$. Associativity of twisted convolution then gives the stated identity.

\textit{(ii)} This is the identity 
\begin{align*}
\varphi^{n-1}_{i,\lambda} *_{-\lambda} \varphi^{n-1}_{j,\lambda} = \delta_{ij}(2\pi)^n|\lambda|^{-n}\varphi^{n-1}_{j,\lambda}
\end{align*}
from \Cref{sec:laguerre-hermite}, with $i=j=k$.

\textit{(iii)} Writing out the definitions,
\begin{align}
W_{S,\lambda}(\tau_c g)(w) &= \int_{\mathbb{C}^n} \tau_c g(w-y)\,\varphi^{n-1}_{k,\lambda}(y)\,
e^{-i\frac{\lambda}{2}\mathrm{Im}(w\cdot\bar y)}\,dy \notag\\
&= \int_{\mathbb{C}^n} e^{i\frac{\lambda}{2}\mathrm{Im}((w-y)\cdot\bar c)} g(w-y+c)\,
\varphi^{n-1}_{k,\lambda}(y)\, e^{-i\frac{\lambda}{2}\mathrm{Im}(w\cdot\bar y)}\, dy \notag\\
&= e^{i\frac{\lambda}{2}\mathrm{Im}(w\cdot\bar c)}\,
\int_{\mathbb{C}^n} g(w+c-y)\,\varphi^{n-1}_{k,\lambda}(y)\,
e^{-i\frac{\lambda}{2}\mathrm{Im}(y\cdot\bar c)}\,e^{-i\frac{\lambda}{2}\mathrm{Im}(w\cdot\bar y)}\,dy \notag\\
&=: e^{i\frac{\lambda}{2}\mathrm{Im}(w\cdot\bar c)}\, I(w) \label{eq:lhs}
\end{align}
On the other hand,
\begin{align}
W_{S,\lambda}(M_c g)(w+c) &= \int_{\mathbb{C}^n} e^{i\lambda\,\mathrm{Im}((w+c-y)\cdot\bar c)}
g(w+c-y)\,\varphi^{n-1}_{k,\lambda}(y)\, e^{-i\frac{\lambda}{2}\mathrm{Im}((w+c)\cdot\bar y)}\,dy \notag\\
&= e^{i\lambda\,\mathrm{Im}(w\cdot\bar c)}\, I(w) \label{eq:rhs}
\end{align}
From Equations (\ref{eq:lhs}) and (\ref{eq:rhs}), we get
\[
e^{i\frac{\lambda}{2}\mathrm{Im}(c\cdot\bar w)}\, W_{S,\lambda}(M_c g)(w+c)
= e^{i\frac{\lambda}{2}\mathrm{Im}(w\cdot\bar c)}\, I(w) = W_{S,\lambda}(\tau_c g)(w),
\]
which is the stated identity.\\

(iv) Substituting $y'=\sigma y$, which is measure-preserving since $\sigma \in U(n)\subset O(2n)$,
\begin{align*}
W_{S,\lambda}(\rho_\sigma g)(z) &= \int_{\mathbb{C}^n} \rho_\sigma g(z-y)\,\varphi^{n-1}_{k,\lambda}(y)\, e^{-i\frac{\lambda}{2}\mathrm{Im}(z\cdot\bar y)}\,dy \\
&= \int_{\mathbb{C}^n} g(\sigma z - y')\,\varphi^{n-1}_{k,\lambda}(\sigma^{-1}y')\, e^{-i\frac{\lambda}{2}\mathrm{Im}(z\cdot\overline{\sigma^{-1}y'})}\,dy' \\
&= \int_{\mathbb{C}^n} g(\sigma z - y')\,\varphi^{n-1}_{k,\lambda}(y')\, e^{-i\frac{\lambda}{2}\mathrm{Im}(\sigma z\cdot\bar y')}\,dy' \\
&= W_{S,\lambda}(g)(\sigma z),
\end{align*}
where the third equality uses two facts. First, $\varphi^{n-1}_{k,\lambda}$ is radial, hence $U(n)$-invariant, so $\varphi^{n-1}_{k,\lambda}(\sigma^{-1}y') = \varphi^{n-1}_{k,\lambda}(y')$. Second, $\sigma$ is unitary, so $\mathrm{Im}(z\cdot\overline{\sigma^{-1}y'}) = \mathrm{Im}(\sigma z\cdot \bar y')$.
\end{proof}

For $r>0$, let $\delta_r f(z,t) := f(rz, r^2t)$ denote the anisotropic dilation on $\mathbb{H}^n$.

\begin{thm}[Dilation covariance]
\label{thm:dilation-covariance}
For $f \in L^1\cap L^2(\mathbb{H}^n)$, $a = (\lambda,(2k+n)|\lambda|) \in R_k$, and $w \in \mathbb{C}^n$,
\begin{equation*}
\mathcal{F}_S(\delta_r f)(a,w) = r^{-(2n+2)}\, \mathcal{F}_S(f)\big(r^{-2}a,\, rw\big).
\end{equation*}
\end{thm}

\begin{proof}
Substituting $t'=r^2t$ in
\begin{align*}
(\delta_r f)^{-\lambda}(z) = \int_{\mathbb{R}} e^{-i\lambda t} f(rz,r^2t)\,dt
\end{align*}
gives
\begin{align}
(\delta_r f)^{-\lambda}(z) = r^{-2} f^{-\lambda/r^2}(rz). \label{eq:dilated-fibre}
\end{align}
By the definition of $\mathcal{F}_S(\delta_r f)(a,w)$ as a twisted convolution and \eqref{eq:dilated-fibre},
\begin{align}
\mathcal{F}_S(\delta_r f)(a,w)
= r^{-2}\int_{\mathbb{C}^n} f^{-\lambda/r^2}(rw-ry)\,\varphi^{n-1}_{k,\lambda}(y)\,
e^{-i\frac{\lambda}{2}\mathrm{Im}(w\cdot\bar y)}\, dy. \label{eq:after-fibre-sub}
\end{align}
A change of variable and an application of the scaling identities
\begin{align*}
\varphi^{n-1}_{k,\lambda}(y'/r) = \varphi^{n-1}_{k,\lambda/r^2}(y'), \qquad
\mathrm{Im}\big(w\cdot\overline{y'/r}\big) = \frac{1}{r^2}\,\mathrm{Im}(rw\cdot\bar y'),
\end{align*}
turns \eqref{eq:after-fibre-sub} into
\begin{align}
\mathcal{F}_S(\delta_r f)(a,w)
&= r^{-(2n+2)} \int_{\mathbb{C}^n} f^{-\lambda/r^2}(rw-y')\,\varphi^{n-1}_{k,\lambda/r^2}(y')\,
e^{-i\frac{\lambda/r^2}{2}\mathrm{Im}(rw\cdot\bar y')}\, dy' \notag\\
&= r^{-(2n+2)}\, \big(f^{-\lambda/r^2} *_{-\lambda/r^2} \varphi^{n-1}_{k,\lambda/r^2}\big)(rw). \label{eq:final-conv}
\end{align}
Writing $\lambda' := r^{-2}\lambda$, so that $-\lambda/r^2 = -\lambda'$ and $a' := (\lambda',(2k+n)|\lambda'|) = r^{-2}a$, the right-hand side of \eqref{eq:final-conv} is $\mathcal{F}_S(f)(r^{-2}a, rw)$, which is the stated identity.
\end{proof}

\begin{rmk}
Dilation covariance is a property of $\mathcal{F}_S$, realised through the joint dilation $\delta_r(z,t)=(rz,r^2t)$ on $\mathbb{H}^n$. A dilation $z\mapsto rz$ on $\mathbb{C}^n$ alone produces a phase mismatch in twisted convolution, so it is not a property of $W_{S,\lambda}$ in isolation. Unitary invariance, by contrast, is already a property of $W_{S,\lambda}$ at the $\mathbb{C}^n$ level, requiring no lift to $\mathbb{H}^n$.
\end{rmk}

\subsection{Idempotency, self-adjointness, and the Plancherel bound}

\begin{prop}
\label{prop:idempotent-self-adjoint}
\begin{enumerate}
\item[(a)] $W_{S,\lambda} \circ W_{S,\lambda} = (2\pi)^n|\lambda|^{-n}\, W_{S,\lambda}$.
\item[(b)] $(W_{S,\lambda}u,v) = (u,W_{S,\lambda}v)$ for $u,v \in L^2(\mathbb{C}^n)$.
\item[(c)] For $u \in L^2(\mathbb{C}^n)$,
\begin{equation} \label{eq:swtbound}
\|W_{S,\lambda}(u)\|_2^2 = (2\pi)^{2n}|\lambda|^{-2n}\, \|P_ku\|_2^2
\le (2\pi)^{2n}|\lambda|^{-2n}\, \|u\|_2^2,
\end{equation}
with equality if and only if $u$ lies in the $k$-th eigenspace of $L_{-\lambda}$.
\end{enumerate}
\end{prop}

\begin{proof}
Immediate from $W_{S,\lambda} = (2\pi)^n|\lambda|^{-n}P_k$ and $P_k$
being an orthogonal projection.
\end{proof}

\section{New Results for Strichartz Fourier transform via Strichartz Weyl transform}
\label{sec:new-results}

The results in this section are obtained directly from the
identification $W_{S,\lambda} = (2\pi)^n|\lambda|^{-n}P_k$ of
\Cref{sec:swt}, independently of any further results built from it.
Each has no stated counterpart for $\mathcal{F}_S$ prior to this paper.

\subsection{A sharp $L^p \to L^2$ bound at fixed $k$}
\label{sec:sharp-lp}
\subsubsection{Worked examples}
\label{sec:worked-examples}

We illustrate two simple examples.
\begin{exam}
$W_{S,\lambda}\big(\varphi^{n-1}_{j,\lambda}\big) = \delta_{jk}\,
(2\pi)^n|\lambda|^{-n}\, \varphi^{n-1}_{k,\lambda}$.
\end{exam}

\begin{exam}
Using the Moyal identity with $\mu=\nu=\gamma$ and summing over
$|\gamma| = k$,
\begin{equation*}
W_{S,\lambda}\big(\Phi^\lambda_{\alpha\beta}\big) = \delta_{|\beta|,k}\,
(2\pi)^n|\lambda|^{-n}\, \Phi^\lambda_{\alpha\beta},
\end{equation*}
nonzero for every $\alpha$ once $|\beta| = k$.
\end{exam}

\begin{exam}
For $g(z) = e^{-c|z|^2}$, $c>0$, a direct Gaussian computation gives, at
$k=0$,
\begin{equation*}
W_{S,\lambda}\big(e^{-c|\cdot|^2}\big)(z) =
\left(\frac{4\pi}{4c+|\lambda|}\right)^n \varphi^{n-1}_{0,\lambda}(z),
\end{equation*}
recovering $(2\pi)^n|\lambda|^{-n}$ at $c = |\lambda|/4$, consistent
with the preceding example. The computation completes the square in
the twisted convolution integral, and the cross terms cancel exactly.
\end{exam}

\begin{thm}
\label{thm:sharp-lp-l2}
Let $1 \le p \le \frac{2n}{n+1}$, $\lambda \in \mathbb{R}^*$, $k \in
\mathbb{N}$. For $f \in L^p(\mathbb{C}^n)$,
\begin{equation*}
\|W_{S,\lambda,k}(f)\|_2 \le C\,(2\pi)^n\, |\lambda|^{\frac{n}{p}-\frac{3n}{2}}\,
k^{n(\frac{1}{p}-\frac{1}{2})-\frac{1}{2}}\, \|f\|_p,
\end{equation*}
with $C = C(n,p)$ independent of $k$ and $\lambda$. The exponent of $k$
is sharp.
\end{thm}

\begin{proof}
\emph{Step 1 ($\lambda = 1$).}
Write $\varphi^{n-1}_k(z) := (2\pi)^{-n}\varphi^{n-1}_{k,1}(z)$ and set
$Q_k f := \varphi^{n-1}_k *_1 f$.
By \cite[Proposition 4.2]{ThangaveluWeylMultipliers1991},
\begin{equation}
\label{eq:sharp-cited}
\|Q_k f\|_2 \le C\, k^{n(1/p - 1/2) - 1/2}\, \|f\|_p,
\qquad 1 \le p \le \frac{2n}{n+1}.
\end{equation}
By the identity $f *_\mu g = g *_{-\mu} f$ at $\mu = 1$, and the
definition of $W_{S,1}$,
\begin{equation*}
\label{eq:sharp-identification}
Q_k f = \varphi^{n-1}_k *_1 f = f *_{-1} \varphi^{n-1}_k
= (2\pi)^{-n} \big( f *_{-1} \varphi^{n-1}_{k,1} \big)
= (2\pi)^{-n} W_{S,1}(f).
\end{equation*}
This gives $W_{S,1,k}(f) = (2\pi)^n Q_k f$.
Substituting into \eqref{eq:sharp-cited},
\begin{equation}
\label{eq:sharp-step1}
\|W_{S,1,k}(f)\|_2 \le C\,(2\pi)^n\, k^{n(1/p-1/2)-1/2}\, \|f\|_p,
\qquad 1 \le p \le \frac{2n}{n+1}.
\end{equation}

\emph{Step 2 ($\lambda = -1$).}
By the conjugation proposition (Proposition~\ref{prop:conjugation}),
with $g = \bar f$,
\[
W_{S,-1}(f) = \overline{W_{S,1}(\bar f)}.
\]
Combined with \eqref{eq:sharp-step1}, this gives,
\begin{align*}
\|W_{S,-1,k}(f)\|_2
&= \|W_{S,1,k}(\bar f)\|_2 \\
&\le C\,(2\pi)^n\, k^{n(1/p-1/2)-1/2}\, \|\bar f\|_p \\
&= C\,(2\pi)^n\, k^{n(1/p-1/2)-1/2}\, \|f\|_p.
\end{align*}
Thus
\eqref{eq:sharp-step1} holds verbatim with $1$ replaced by $-1$.\\

\emph{Step 3 (scaling to general $\lambda$).}
From the defining formula, $\varphi^{n-1}_{k,\lambda}(z) =
\varphi^{n-1}_{k,\operatorname{sgn}(\lambda)}(\sqrt{|\lambda|}z)$.
Substituting $w = \sqrt{|\lambda|}w'$ in the definition of $W_{S,\lambda}$,
and rescaling the twisted phase accordingly, gives, writing
$\delta_r f(z) := f(rz)$,
\[
W_{S,\lambda}(f)(z) = |\lambda|^{-n}\, W_{S,\operatorname{sgn}(\lambda),k}
\big(\delta_{1/\sqrt{|\lambda|}} f\big)\big(\sqrt{|\lambda|}z\big).
\]
This implies,
\begin{align*}
\|W_{S,\lambda,k}(f)\|_2
&= |\lambda|^{-n}\cdot|\lambda|^{-n/2}\,
\big\|W_{S,\operatorname{sgn}(\lambda),k}\big(\delta_{1/\sqrt{|\lambda|}}f\big)\big\|_2 \\
&= |\lambda|^{-3n/2}\,
\big\|W_{S,\operatorname{sgn}(\lambda),k}\big(\delta_{1/\sqrt{|\lambda|}}f\big)\big\|_2.
\end{align*}
Applying \eqref{eq:sharp-step1}, or its $\lambda=-1$ counterpart, to
$\delta_{1/\sqrt{|\lambda|}}f$, and using
$\|\delta_{1/\sqrt{|\lambda|}}f\|_p = |\lambda|^{n/p}\|f\|_p$,
\begin{align*}
\|W_{S,\lambda,k}(f)\|_2 &\le |\lambda|^{-3n/2}\cdot C\,(2\pi)^n\,
k^{n(1/p-1/2)-1/2} \cdot |\lambda|^{n/p}\, \|f\|_p \\
&= C\,(2\pi)^n\, |\lambda|^{n/p-3n/2}\, k^{n(1/p-1/2)-1/2}\, \|f\|_p,
\end{align*}
which is the stated bound.\\

\emph{Step 4 (sharpness).}
We first identify the relevant object in \cite{KochRicci2007}. Write
\[
L^{\mathrm{tw}} := -\tfrac12\sum_{j=1}^n\big[(\partial_{x_j}+iy_j)^2+(\partial_{y_j}-ix_j)^2\big]
\]
for the twisted Laplacian of \cite{KochRicci2007}. Expanding each square
and summing over $j$ gives
\[
L^{\mathrm{tw}} = -\tfrac12\Delta + \tfrac12|z|^2 + iN,
\]
where $N$ is the rotation operator of Section~\ref{sec:laguerre-hermite}.
Comparing with $L_\mu = -\Delta+\tfrac14\mu^2|z|^2-i\mu N$ at $\mu=2$,
\[
\tfrac12 L_2 = -\tfrac12\Delta+\tfrac12|z|^2 - iN,
\]
which agrees with $L^{\mathrm{tw}}$ except in the sign of the
$N$-term. Since $N$ annihilates radial functions, $L^{\mathrm{tw}}u_0 =
\tfrac12 L_2 u_0$ for every radial $u_0$. On eigenvalues,
$L^{\mathrm{tw}}$ has spectrum $n+2k$ and $\tfrac12L_2$ has spectrum
$\tfrac12\cdot 2(2k+n) = n+2k$, so the two coincide. Hence the radial
eigenfunction $f_k$ of \cite[Section 2]{KochRicci2007}, that is, the unique (up
to scalar) radial eigenfunction of $L^{\mathrm{tw}}$ with eigenvalue
$n+2k$, coincides with the unique radial eigenfunction of $L_{-2}$ at
the same eigenvalue, namely $\varphi^{n-1}_{k,2}$.

By the first worked example of Section~\ref{sec:worked-examples}
applied at $\lambda=2$,
\[
W_{S,2,k}(\varphi^{n-1}_{k,2}) = (2\pi)^n 2^{-n}\varphi^{n-1}_{k,2}.
\]
Testing the general-$\lambda$ bound established at the end of Step 3,
at $\lambda=2$, against $f = \varphi^{n-1}_{k,2}$, reduces sharpness of
the $k$-exponent to the asymptotic ratio
$\|\varphi^{n-1}_{k,2}\|_{p'}/\|\varphi^{n-1}_{k,2}\|_2$ as
$k\to\infty$. Equivalently, this is the ratio $\|f_k\|_{p'}/\|f_k\|_2$
for this same extremal.

Write $\varrho(q) = \tfrac{d-2}2-\tfrac dq$, $d=2n$. The sharp
far-range exponent of \cite{KochRicci2007} holds for $q \ge
\tfrac{2(d+1)}{d-1}$. Dualising via $\|P_k\|_{L^p\to L^2} =
\|P_k\|_{L^2\to L^{p'}}$, with $q=p'$ and $\tfrac1{p'}=1-\tfrac1p$,
gives
\[
\|f_k\|_{p'}/\|f_k\|_2 \ \gtrsim\ k^{\varrho(p')/2}
= k^{\,n/p - (n+1)/2} = k^{\,n(1/p-1/2)-1/2},
\]
which is the $k$-exponent of the theorem.

This far range $q \ge \tfrac{2(d+1)}{d-1}$ translates, via $q=p'$, to
$p \le \tfrac{4n+2}{2n+3}$. This range contains the range $p \le
\tfrac{2n}{n+1}$ assumed above, since $(n+1)(4n+2) - 2n(2n+3) = 2 > 0$
for every $n$. This shows the $k$-exponent is attained, up to
constants, at $\lambda=2$.

Sharpness at every other nonzero $\lambda$, including the general
$\lambda$ of the theorem, then follows by the same dilation identity
of Step 3, which is an exact equality and so transfers sharpness, not
merely boundedness.
\end{proof}

\subsection{Fock-Bargmann transport}
\label{sec:fock-bargmann}

\begin{defn}
The Fock space $\mathcal{F}_n$ is the space of entire functions
$F : \mathbb{C}^n \to \mathbb{C}$ satisfying
\[
\|F\|^2_{\mathcal{F}_n} := \pi^{-n} \int_{\mathbb{C}^n} |F(\zeta)|^2
e^{-|\zeta|^2}\, d\zeta < \infty,
\]
equipped with this norm and the associated inner product.
\end{defn}

For $\alpha \in \mathbb{N}^n$, write
\[
\Phi^\lambda_\alpha(\xi) := |\lambda|^{n/4}\Phi_\alpha(\sqrt{|\lambda|}\xi),
\qquad
e_\alpha(\zeta) := \frac{\zeta^\alpha}{\sqrt{\alpha!}}.
\]
Let $B_\lambda : L^2(\mathbb{R}^n) \to \mathcal{F}_n$ be the
$\lambda$-adapted Bargmann transform \cite{Bargmann1961}, the unitary
isomorphism determined by
\[
B_\lambda \Phi^\lambda_\alpha = e_\alpha, \qquad \alpha \in \mathbb{N}^n.
\]
Let
\[
Q^{(\lambda)}_k := \text{orthogonal projection of } L^2(\mathbb{R}^n)
\text{ onto } \operatorname{span}\{\Phi^\lambda_\alpha : |\alpha|=k\},
\]
\[
\mathcal{F}_n^{(k)} := \operatorname{span}\{e_\alpha : |\alpha|=k\}
\quad (\text{the degree-$k$ homogeneous polynomials}),
\]
\[
\Pi_k := \text{orthogonal projection of } \mathcal{F}_n \text{ onto }
\mathcal{F}_n^{(k)}.
\]
Write
\[
\widetilde W_{-\lambda}(f) := B_\lambda\, W_{-\lambda}(f)\, B_\lambda^{-1}.
\]

\begin{lem}[Weyl transform of the Laguerre kernel]\label{lem:weyl-laguerre}
For $\lambda\in\mathbb R^*$ and $k\in\mathbb N$,
\[
W_{-\lambda}(\varphi^{n-1}_{k,\lambda}) = (2\pi)^n|\lambda|^{-n}\,Q_k^{(\lambda)},
\]
where $Q_k^{(\lambda)}$ is the orthogonal projection of $L^2(\mathbb R^n)$ onto
$\mathrm{span}\{\Phi^\lambda_\alpha : |\alpha|=k\}$.
\end{lem}

\begin{proof}
This is the standard connection between the Hermite projections and the
Laguerre functions via the Weyl transform (\cite{ThangaveluLectures1993}).
\end{proof}

\begin{thm}
For $g \in L^1(\mathbb{C}^n)$,
\begin{enumerate}
\item[(i)] $\Pi_k = B_\lambda Q^{(\lambda)}_k B_\lambda^{-1}$.
\item[(ii)] $\widetilde W_{-\lambda}\big(W_{S,\lambda}(g)\big) =
(2\pi)^n|\lambda|^{-n}\, \widetilde W_{-\lambda}(g) \circ \Pi_k$.
\end{enumerate}
\end{thm}

\begin{proof}
By Lemma~\ref{lem:weyl-laguerre},
\[
  W_{-\lambda}(\varphi^{n-1}_{k,\lambda}) = (2\pi)^n|\lambda|^{-n}\, Q^{(\lambda)}_k.
\]

\textit{(i)} Since $B_\lambda$ is unitary, conjugating $Q_k^{(\lambda)}$ by $B_\lambda$ gives the
orthogonal projection on $\mathcal{F}_n$ whose range is
$B_\lambda\big(\operatorname{span}\{\Phi_\alpha^\lambda : |\alpha|=k\}\big)$. This image is
\[
  B_\lambda\big(\operatorname{span}\{\Phi_\alpha^\lambda : |\alpha|=k\}\big)
  = \operatorname{span}\{e_\alpha : |\alpha|=k\}
  = \mathcal{F}_n^{(k)},
\]
which is $\Pi_k$ by definition.

\textit{(ii)} The Weyl transform satisfies the homomorphism property
$W_\mu(f \ast_\mu h) = W_\mu(f) W_\mu(h)$. Taking $\mu = -\lambda$ and recalling
$W_{S,\lambda}(g) = g \ast_{-\lambda} \varphi_{k,\lambda}^{n-1}$, this gives
\[
  W_{-\lambda}\big(W_{S,\lambda}(g)\big) = W_{-\lambda}(g) \cdot W_{-\lambda}(\varphi_{k,\lambda}^{n-1}).
\]
By Lemma~\ref{lem:weyl-laguerre}, the right-hand side equals
$(2\pi)^n|\lambda|^{-n}\, W_{-\lambda}(g)\, Q_k^{(\lambda)}$. Conjugating both sides by $B_\lambda$
and using part (i), this becomes
\[
  \widetilde W_{-\lambda}\big(W_{S,\lambda}(g)\big)
  = (2\pi)^n|\lambda|^{-n}\, \widetilde W_{-\lambda}(g) \circ \Pi_k,
\]
which is the stated identity.
\end{proof}

\begin{thm}
For $F \in \mathcal{F}_n$,
\begin{enumerate}
\item[(i)] $\Pi_k^2 = \Pi_k = \Pi_k^*$ on $\mathcal{F}_n$.
\item[(ii)] $\|\Pi_k F\|_{\mathcal{F}_n} \le \|F\|_{\mathcal{F}_n}$, with equality if and only if
$F \in \mathcal{F}_n^{(k)}$.
\item[(iii)]
\begin{equation*}
F = \sum_{k=0}^\infty \Pi_k F \quad \text{in } \mathcal{F}_n, \qquad \text{and} \qquad
\|F\|^2_{\mathcal{F}_n} = \sum_{k=0}^\infty \|\Pi_k F\|^2_{\mathcal{F}_n}.
\end{equation*}
\end{enumerate}
\end{thm}

\begin{proof}
Since $B_\lambda$ is unitary, it transports every identity already established for $P_k$ onto
the corresponding identity for $\Pi_k$, using $\Pi_k = B_\lambda Q^{(\lambda)}_k B_\lambda^{-1}$.

\textit{(i)} Idempotency and self-adjointness of $P_k$ carry over directly, giving
$\Pi_k^2 = \Pi_k = \Pi_k^*$.

\textit{(ii)} The Plancherel bound for $P_k$ transports to
$\|\Pi_k F\|_{\mathcal{F}_n} \le \|F\|_{\mathcal{F}_n}$, with equality if and only if
$F \in \mathcal{F}_n^{(k)}$.

\textit{(iii)} Completeness of the special Hermite expansion shows that
$F = \sum_{k=0}^\infty \Pi_k F$ in $\mathcal{F}_n$. Summing the bound of part (ii) over $k$
then establishes $\|F\|^2_{\mathcal{F}_n} = \sum_{k=0}^\infty \|\Pi_k F\|^2_{\mathcal{F}_n}$.
\end{proof}

\begin{rmk}
This Fock-space transport has no stated counterpart for
$\mathcal{F}_S$ in \cite{Thangavelu2023}, which does not treat the
Fock-space realisation of $\pi_\lambda$. Whether a corresponding
fan-level statement for $\mathcal{F}_S$ holds is accordingly left
pending, and is not addressed in the present paper. We record this
Fock-space transport as a structural result about $W_{S,\lambda}$ in
its own right, in the same spirit as the Fock-space realisation of the
classical Weyl transform, which is foundational in its own literature
independently of any single downstream use. It is not applied
elsewhere in the present paper.
\end{rmk}

\subsection{A Benedicks-Type Uncertainty Theorem}
\label{sec:benedicks}

As an application of the identification of $W_{S,\lambda}$ with a
spectral projection of the elliptic operator $L_{-\lambda}$, we record
a qualitative uncertainty principle in the spirit of Benedicks'
classical theorem for the Euclidean Fourier transform.

\begin{lem}[Analytic hypoellipticity, \cite{MorreyNirenberg1957}]
Let $L$ be an elliptic differential operator on $\mathbb{R}^N$ with
real-analytic coefficients. If $u$ is a distributional solution of $Lu
= cu$ for a constant $c$, then $u$ is real-analytic on $\mathbb{R}^N$.
\end{lem}

\begin{thm}
\label{thm:benedicks}
Fix $\lambda \in \mathbb{R}^*$. Let $g \in L^2(\mathbb{C}^n)$. If, for
every $k \in \mathbb{N}$, the set $\{z \in \mathbb{C}^n :
W_{S,\lambda,k}(g)(z) \ne 0\}$ has finite Lebesgue measure, then $g
\equiv 0$.
\end{thm}

\begin{proof}
Fix $k \in \mathbb{N}$ and write $g_k := P_k g$. Then
\[
  W_{S,\lambda,k}(g) = (2\pi)^n|\lambda|^{-n} g_k.
\]
Since $g_k$ lies in the $k$-th eigenspace of $L_{-\lambda}$,
\[
  L_{-\lambda}\big(W_{S,\lambda,k}(g)\big) = (2k+n)|\lambda|\, W_{S,\lambda,k}(g).
\]
The operator $L_{-\lambda}$ is elliptic on $\mathbb{C}^n \cong \mathbb{R}^{2n}$ with polynomial,
hence real-analytic, coefficients. So by the preceding lemma $W_{S,\lambda,k}(g)$ is
real-analytic on $\mathbb{R}^{2n}$.

Suppose $W_{S,\lambda,k}(g) \not\equiv 0$. Let
\[
  Z := \{z \in \mathbb{C}^n : W_{S,\lambda,k}(g)(z) = 0\}.
\]
Since $W_{S,\lambda,k}(g)$ is real-analytic and not identically zero, $Z$ has Lebesgue
measure zero. Hence the measure of its complement is not finite, contradicting the
hypothesis that $\{z \in \mathbb{C}^n : W_{S,\lambda,k}(g)(z) \ne 0\}$ has finite measure.

Thus $W_{S,\lambda,k}(g) \equiv 0$, that is, $g_k = P_k g = 0$. Since $k \in \mathbb{N}$ was
arbitrary, $P_k g = 0$ for every $k$, and completeness of the special Hermite expansion gives
\[
  g = \sum_{k=0}^\infty P_k g = 0.
\]
\end{proof}

\begin{rmk}
Real-analyticity of $W_{S,\lambda,k}(g)$ at each $k$ follows automatically from its lying in
the fixed eigenspace of the elliptic operator $L_{-\lambda}$, regardless of any support
condition on $g$. This is a structural difference from the classical Euclidean Benedicks
theorem, where real-analyticity of $\widehat f$ is not automatic and is instead a consequence
of assuming $f$ compactly supported.
\end{rmk}

\begin{rmk}
Theorem~\ref{thm:benedicks} fixes $\lambda$ and lets only $k$ vary. A genuinely broader
Benedicks-type theorem holds for $\mathcal{F}_S$ itself, with a hypothesis imposed jointly
across $\lambda$ and $k$ simultaneously. This stronger statement, due to Dabra, Dasgupta, and
Gulia \cite{DabraDasguptaGulia2026}, lies beyond what the fixed-ray method can reach. We
discuss this in detail in \Cref{sec:limits}.
\end{rmk}

\section{Core Toolkit: Range, Inversion, Plancherel, and Continuity}
\label{sec:core-toolkit}

This section establishes two layers of results. First, fixed-$\lambda$ results about $W_{S,\lambda}$ itself, building on the fixed-ray results of \Cref{sec:swt}. These include a Plancherel identity, an inversion formula, a self-reproducing range condition, a Hausdorff-Young inequality, and a Riemann-Lebesgue/continuity theory for the Strichartz Weyl transform. These results are new, since $W_{S,\lambda}$ is introduced in this paper and has no prior literature. Second, assembling these fixed-ray results across all of $\Omega$ recovers the corresponding toolkit for the Strichartz Fourier transform itself, namely an inversion formula, a Plancherel theorem, a characterisation of its range, a Hausdorff-Young inequality, and a Riemann-Lebesgue/continuity theory. This $\mathcal{F}_S$-level toolkit is not new. It was established directly, by Laguerre and Bessel asymptotics, in \cite{Thangavelu2023}. We include it here because idempotency, self-adjointness, and the Plancherel bound for a single orthogonal projection already supply, as formal consequences, a body of theorems that were originally established in \cite{Thangavelu2023} through direct estimation involving Laguerre and Bessel asymptotics.

\subsection{Range, inversion, and Plancherel}

\subsubsection{The self-reproducing range condition}

\begin{thm}
\label{thm:self-reproducing}
A function $h \in \mathcal{S}(\mathbb{C}^n)$ lies in the $k$-th eigenspace of $L_{-\lambda}$ if and only if
\begin{equation*}
h *_{-\lambda} \varphi^{n-1}_{k,\lambda} = (2\pi)^n|\lambda|^{-n}\, h.
\end{equation*}
In particular, $W_{S,\lambda}$ maps $\mathcal{S}(\mathbb{C}^n)$ into itself, and $W_{S,\lambda}(g)$ satisfies this condition for every $g \in \mathcal{S}(\mathbb{C}^n)$.
\end{thm}

\begin{proof}
The stated equivalence is $P_kh=h$, rewritten via \eqref{eq:laguerre-projection}. Twisted convolution against a Schwartz kernel preserves $\mathcal{S}(\mathbb{C}^n)$, so $W_{S,\lambda}$ maps $\mathcal{S}(\mathbb{C}^n)$ into itself. Also, $W_{S,\lambda}(g)$ satisfies this condition for every $g\in\mathcal{S}(\mathbb{C}^n)$, by the idempotency relation $W_{S,\lambda}\circ W_{S,\lambda}=(2\pi)^n|\lambda|^{-n}W_{S,\lambda}$ applied to $g$.
\end{proof}

\begin{rmk}
The order of the kernel matters in this condition, and it must appear on the right of $*_{-\lambda}$. Write $h = \Phi^\lambda_{\rho\sigma}$ with $|\sigma|=k$, so that $h$ lies in the $k$-th eigenspace. The Moyal identity gives
\begin{align*}
h *_{-\lambda} \varphi^{n-1}_{k,\lambda} = (2\pi)^n|\lambda|^{-n}h.
\end{align*}
On the other hand,
\begin{align*}
\varphi^{n-1}_{k,\lambda} *_{-\lambda} h = \delta_{|\rho|,k}\,(2\pi)^n|\lambda|^{-n}h,
\end{align*}
which vanishes unless $|\rho|=k$ as well. This is not satisfied by a general element of the eigenspace, since $h$ lying in the $k$-th eigenspace only constrains $\sigma$, not $\rho$.
\end{rmk}

\subsubsection{Plancherel identity and inversion at fixed $\lambda$}

\begin{thm}
\label{thm:plancherel-fixed-lambda}
Fix $\lambda \in \mathbb{R}^*$. For every $u \in L^2(\mathbb{C}^n)$,
\begin{equation*}
\|u\|_2^2 = (2\pi)^{-2n}|\lambda|^{2n} \sum_{k=0}^\infty \|W_{S,\lambda,k}(u)\|_2^2.
\end{equation*}
\end{thm}

\begin{proof}
By \eqref{eq:swtbound},
\begin{align*}
\|W_{S,\lambda,k}(u)\|_2^2 = (2\pi)^{2n}|\lambda|^{-2n}\|P_ku\|_2^2.
\end{align*}
Summing over $k$ and using Parseval's identity for the special Hermite expansion, $\sum_k \|P_ku\|_2^2 = \|u\|_2^2$, gives the stated identity.
\end{proof}

\begin{cor}[Inversion at fixed $\lambda$]
\label{cor:inversion-fixed-lambda}
For $f \in L^1\cap L^2(\mathbb{H}^n)$,
\begin{equation*}
f(z,t) = (2\pi)^{-n-1}\int_{\mathbb{R}} e^{i\lambda t}\,|\lambda|^n \sum_{k=0}^\infty \mathcal{F}_S(f)(a_k,z)\, d\lambda, \qquad a_k = (\lambda,(2k+n)|\lambda|),
\end{equation*}
in the $L^2$ sense.
\end{cor}

\begin{proof}
Fix $\lambda \in \mathbb{R}^*$. By \Cref{thm:plancherel-fixed-lambda},
\begin{equation*}
u = (2\pi)^{-n}|\lambda|^n \sum_{k=0}^\infty W_{S,\lambda,k}(u) \qquad \text{in } L^2(\mathbb{C}^n),
\end{equation*}
for every $u \in L^2(\mathbb{C}^n)$. Taking $u = f^{-\lambda}$ gives
\begin{equation*}
f^{-\lambda}(z) = (2\pi)^{-n}|\lambda|^n \sum_{k=0}^\infty \mathcal{F}_S(f)(a_k,z) \qquad \text{in } L^2(\mathbb{C}^n).
\end{equation*}
Recovering $f(z,t)$ from $f^{-\lambda}$ by Fourier inversion in $t$,
\begin{equation*}
f(z,t) = \frac{1}{2\pi}\int_{\mathbb{R}} e^{i\lambda t} f^{-\lambda}(z)\, d\lambda,
\end{equation*}
and substituting the preceding identity gives the stated formula.
\end{proof}

\subsubsection{The inversion formula and Plancherel theorem over $\Omega$}
The inversion formula for $\mathcal{F}_S$ is due to \cite[Thm.~7.2, first part]{Thangavelu2023}. There, it is obtained by direct computation with the special Hermite expansion at each fixed $\lambda$, verified coefficient-by-coefficient in $k$. Below we reach it instead by assembling the fixed-$\lambda$ inversion identity of \Cref{cor:inversion-fixed-lambda} across $\Omega$, itself a consequence of Parseval's identity for the special Hermite expansion packaged through the operator identity $W_{S,\lambda} = (2\pi)^n|\lambda|^{-n}P_k$.

\begin{thm}[Inversion formula]
\label{thm:inversion-omega}
For any Schwartz class function $f$ on $\mathbb{H}^n$,
\begin{equation*}
f(z,t) = \int_\Omega \int_{\mathbb{C}^n} \widehat f(a,w)\, e_a\big((-w,0)(z,t)\big)\, dw\, d\nu(a).
\end{equation*}
\end{thm}

\begin{proof}
By Corollary~\ref{cor:inversion-fixed-lambda}, for almost every
$\lambda$,
\begin{equation*}
f^{-\lambda}(z) = (2\pi)^{-n}|\lambda|^n \sum_{k=0}^\infty
\widehat f(a,z), \qquad a = (\lambda,(2k+n)|\lambda|) \in R_k,
\end{equation*}
with the sum convergent in $L^2(\mathbb{C}^n)$. Recovering $f(z,t)$
from $f^{-\lambda}$ by Fourier inversion in $t$,
\begin{align}
f(z,t) &= \frac{1}{2\pi}\int_{-\infty}^\infty e^{i\lambda t} f^{-\lambda}(z)\, d\lambda \notag\\
&= (2\pi)^{-n-1}\int_{-\infty}^\infty e^{i\lambda t}\Big(\sum_{k=0}^\infty \widehat f(a,z)\Big)\, |\lambda|^n\, d\lambda. \label{eq:f-from-fhat}
\end{align}
Recalling $e_a(z,t) = e^{i\lambda t}\varphi^{n-1}_{k,\lambda}(z)$ 
we get
\begin{align*}
e_a\big((-w,0)(z,t)\big) = e^{i\lambda t}\, \varphi^{n-1}_{k,\lambda}(z-w)\, e^{-i\frac{\lambda}{2}\mathrm{Im}(w\cdot\bar z)}.
\end{align*}
This gives
\begin{align*}
\int_{\mathbb{C}^n} \widehat f(a,w)\, e_a\big((-w,0)(z,t)\big)\, dw
&= e^{i\lambda t}\, \big(\widehat f(a,\cdot) *_{-\lambda} \varphi^{n-1}_{k,\lambda}\big)(z) \\
&= e^{i\lambda t}\, (2\pi)^n|\lambda|^{-n}\, \widehat f(a,z),
\end{align*}
where the last equality is \Cref{thm:self-reproducing}, since $\widehat f(a,\cdot)$ lies in the $k$-th eigenspace of $L_{-\lambda}$. Integrating over $\Omega$ against $d\nu(a)$,
\begin{align*}
\int_\Omega \int_{\mathbb{C}^n} \widehat f(a,w)\,e_a\big((-w,0)(z,t)\big)\,dw\,d\nu(a)
&= \int_\Omega e^{i\lambda t}\,(2\pi)^n|\lambda|^{-n}\, \widehat f(a,z)\,d\nu(a) \\
&= (2\pi)^{-n-1}\int_{-\infty}^\infty e^{i\lambda t}\,|\lambda|^n \sum_{k=0}^\infty \widehat f(a,z)\, d\lambda,
\end{align*}
which is the right-hand side of \Cref{eq:f-from-fhat}.
\end{proof}

The Plancherel formula for $\mathcal{F}_S$ is established in \cite[Thm.~7.2, second part]{Thangavelu2023}, where it is obtained from the Plancherel identity for the special Hermite expansion combined with the classical Plancherel theorem for the Euclidean Fourier transform in $t$. The proof given here follows the same two-step combination. The only difference is that the special Hermite Plancherel identity used there is replaced by the fixed-$\lambda$ identity of Theorem~\ref{thm:plancherel-fixed-lambda}, which is itself a restatement of the Plancherel bound of a single orthogonal projection, summed over $k$.

\begin{thm}[Plancherel formula]
\label{thm:plancherel-omega}
For $f \in L^1 \cap L^2(\mathbb{H}^n)$,
\begin{equation} \label{eq:sftplancherel}
\int_{\mathbb{H}^n} |f(z,t)|^2\, dz\, dt = \int_\Omega \int_{\mathbb{C}^n}
|\widehat f(a,w)|^2\, dw\, d\nu(a).
\end{equation}
\end{thm}

\begin{proof}
By Theorem~\ref{thm:plancherel-fixed-lambda},
\begin{equation*}
\|f^{-\lambda}\|_{L^2(\mathbb{C}^n)}^2 = (2\pi)^{-2n}|\lambda|^{2n}
\sum_{k=0}^\infty \|\widehat f(a,\cdot)\|_{L^2(\mathbb{C}^n)}^2,
\qquad a = (\lambda,(2k+n)|\lambda|) \in R_k.
\end{equation*}
Combining this with the classical Plancherel theorem for the Euclidean
Fourier transform in $t$,
\begin{equation*}
\int_{\mathbb{H}^n}|f(z,t)|^2\,dz\,dt = \frac{1}{2\pi}\int_{-\infty}^\infty
\|f^{-\lambda}\|_{L^2(\mathbb{C}^n)}^2\, d\lambda,
\end{equation*}
and recalling the definition of $d\nu$, gives the stated identity.
\end{proof}

\subsubsection{Unitary range of the Strichartz Fourier transform}

The unitary-range characterisation of $\mathcal{F}_S$ is due to \cite[Thm.~7.3]{Thangavelu2023}, where it is verified using the orthogonality of the Laguerre kernel across eigenspaces. Framed through $W_{S,\lambda}$, the same characterisation follows directly from the idempotency and self-adjointness of $P_k$, already established in \Cref{sec:swt} as general facts about spectral projections.

\begin{thm}[Unitary range]
\label{thm:unitary-range}
The Strichartz Fourier transform, initially defined on $L^1\cap
L^2(\mathbb{H}^n)$, extends to a unitary operator from
$L^2(\mathbb{H}^n)$ onto $L^2_0(\Omega\times\mathbb{C}^n,\, d\nu\, dw)$.
\end{thm}

\begin{proof}
Since $P_k$ is an orthogonal
projection, for $h \in L^2(\mathbb{C}^n)$, we have
\[
P_k h = h \quad\Longleftrightarrow\quad h \in \operatorname{Ran}(P_k).
\]
Since $W_{S,\lambda,k} = (2\pi)^n|\lambda|^{-n}P_k$, the fixed-ray
necessary condition at $a \in R_k$, namely
$W_{S,\lambda,k}h = (2\pi)^n|\lambda|^{-n}h$, is just $P_k h = h$.
Hence it identifies membership in $\operatorname{Ran}(W_{S,\lambda,k})$
with membership in $\operatorname{Ran}(P_k)$.

Now let $F \in L^2_0(\Omega\times\mathbb{C}^n, d\nu\, dw)$. As in the
proof of Theorem~\ref{thm:inversion-omega}, define
\[
f^{-\lambda}(z) := (2\pi)^{-n}|\lambda|^n \sum_{k=0}^\infty F(a,z),
\qquad a = (\lambda,(2k+n)|\lambda|),
\]
and recover $f$ from $f^{-\lambda}$ by Fourier inversion in $t$. Applying \Cref{eq:sftplancherel} to
$F$, $f^{-\lambda} \in L^2(\mathbb{C}^n)$ for a.e.\ $\lambda$ and
$f \in L^2(\mathbb{H}^n)$.

Since $F(a,\cdot)$ satisfies the necessary condition at every
$a \in R_k$, the equivalence above gives $P_k F(a,\cdot) = F(a,\cdot)$
for every $k$ and a.e.\ $\lambda$. By \Cref{rmk:comparing-FS}, we get
\[
\widehat f(a,z) = W_{S,\lambda,k}(f^{-\lambda})(z)
= (2\pi)^n|\lambda|^{-n} P_k f^{-\lambda}(z)
= F(a,z).
\]
Thus $\widehat f = F$, completing the proof.
\end{proof}

\begin{rmk}
For $f \in L^2(\mathbb{H}^n)$, $f^\lambda$ is defined only for almost
every $\lambda \in \mathbb{R}^*$, so $\widehat f(a,z)$ may not be
defined at every $a \in \Omega$. In particular, $\widehat f(0,\tau,z)$
need not be defined. As the boundary ray $R_\infty$ has zero Plancherel
measure, this does not affect the preceding theorem.
\end{rmk}

\subsection{The Hausdorff-Young inequality}

\subsubsection{A Hausdorff-Young inequality at fixed $\lambda$}

\begin{thm}
\label{thm:hausdorff-young-fixed}
For $1 \le p \le 2$, $\frac{1}{p} +
\frac{1}{p'} = 1$, there is a measure $\nu^\lambda_p$ on $\mathbb{N}$
such that for every $g \in L^p(\mathbb{C}^n)$,
\begin{equation*}
\left(\sum_{k=0}^\infty \int_{\mathbb{C}^n}
|\widetilde W_{S,\lambda,k}(g)(z)|^{p'}\, dz\, d\nu^\lambda_p(k)\right)^{1/p'}
\le C_{n,p,\lambda}
\left(\int_{\mathbb{C}^n} |g(z)|^p\, dz\right)^{1/p}.
\end{equation*}
\end{thm}

\begin{proof}
At $p=1$: since $|\widetilde\varphi^{n-1}_{k,\lambda}| \le 1$
uniformly in $k$, we have $\|\widetilde W_{S,\lambda,k}(g)\|_\infty \le
\|g\|_1$ for every $k$. Thus $g \mapsto
(\widetilde W_{S,\lambda,k}(g))_k$ is bounded from $L^1(\mathbb{C}^n)$
into $\ell^\infty_k(L^\infty(\mathbb{C}^n))$, with respect to the
counting measure $\nu^\lambda_1$ on $\mathbb{N}$.

\noindent At $p=2$: since $\widetilde W_{S,\lambda,k} = L^{n-1}_k(0)^{-1}
W_{S,\lambda,k}$, the Plancherel bound gives
\begin{equation*}
\|\widetilde W_{S,\lambda,k}(g)\|_2^2 = L^{n-1}_k(0)^{-2}
(2\pi)^{2n}|\lambda|^{-2n}\|P_kg\|_2^2.
\end{equation*}
Summing over $k$ and using completeness,
\begin{equation*}
\sum_{k=0}^\infty L^{n-1}_k(0)^2\, \|\widetilde W_{S,\lambda,k}(g)\|_2^2
= (2\pi)^{2n}|\lambda|^{-2n}\, \|g\|_2^2.
\end{equation*}
That is, $g \mapsto (\widetilde W_{S,\lambda,k}(g))_k$ is, up to the
constant $(2\pi)^n|\lambda|^{-n}$, an isometry from $L^2(\mathbb{C}^n)$
into $\ell^2_k(L^2(\mathbb{C}^n))$, with respect to the weighted measure
$\nu^\lambda_2(k) := L^{n-1}_k(0)^2$.

The measures $\nu_1^\lambda$ and $\nu_2^\lambda$ are mutually
absolutely continuous on the common index set $\mathbb{N}$, so we may
interpolate these two endpoints by Stein-Weiss change-of-measure
interpolation \cite[Thm.~2.9]{SteinWeiss1958}. This gives the stated
inequality for the resulting interpolated measure $\nu_p^\lambda$ on
$\mathbb{N}$.
\end{proof}
\subsubsection{The Hausdorff-Young inequality over $\Omega$}

This is the Hausdorff-Young inequality for $\mathcal{F}_S$
\cite[Thm.~7.4]{Thangavelu2023}. Both there and here, the two
endpoints are combined by the same Stein-Weiss change-of-measure
interpolation \cite{SteinWeiss1958}. What differs is how the two
endpoints are obtained. There, both are established directly against
explicit bounds for the Laguerre kernel and the Bessel function. Here,
both are instead supplied by the uniform bound and the Plancherel
identity already established for $W_{S,\lambda}$. We reprove it by
interpolating between two endpoints drawn from the Strichartz Weyl
transform, namely the uniform bound of the preceding subsection and
the Plancherel identity of Theorem~\ref{thm:plancherel-omega},
carrying out the interpolation itself jointly over $\Omega$, as the
fixed-$\lambda$ interpolation above is confined to a single ray. As a
consequence of the fixed-ray estimate $|\widetilde
W_{S,\lambda,k}(g)(z)| \le \|g\|_{L^1(\mathbb{C}^n)}$, applied with
$g = f^{-\lambda}$ and summed over $k$ and $\lambda$, the Strichartz
Fourier transform of any $f \in L^1(\mathbb{H}^n)$ satisfies
\begin{equation*}
\sup_{(a,z) \in \Omega\times\mathbb{C}^n} |\widetilde f(a,z)| \le c_{n,1}\, \|f\|_1.
\end{equation*}

\begin{thm}[Hausdorff-Young inequality] (\cite{Thangavelu2023})
\label{thm:hausdorff-young-omega}
For $1 \le p \le 2$, there is a measure $\nu_p$ on $\Omega$ such that
for all $f \in L^p(\mathbb{H}^n)$,
\begin{equation*}
\left(\int_\Omega \int_{\mathbb{C}^n} |\widetilde f(a,w)|^{p'}\, dw\,
d\nu_p(a)\right)^{1/p'} \le c_{n,p} \left(\int_{\mathbb{H}^n}
|f(z,t)|^p\, dz\, dt\right)^{1/p}, \qquad \frac{1}{p}+\frac{1}{p'}=1.
\end{equation*}
\end{thm}

\begin{rmk}
The above result follows by assembling Theorem~\ref{thm:hausdorff-young-fixed}
across $\Omega$, in the same manner as the inversion and Plancherel
theorems above. The proof is close in substance to the one given in \cite[Thm.~7.4]{Thangavelu2023}. This involves interpolating
the same two facts via the same Stein-Weiss change-of-measure
interpolation, with the two endpoints now supplied by
$W_{S,\lambda}$'s uniform bound and Plancherel identity rather than
by direct Laguerre and Bessel estimates.
\end{rmk}

\subsection{The Strichartz Weyl transform at the boundary ray}
\label{sec:boundary-ray}

Recall that $R_\infty = \{(0,\tau) : \tau \ge 0\}$ is the boundary of
$\Omega$, where $\mathcal{F}_S$ is given by
\begin{equation*}
\mathcal{F}_S(f)(0,\tau,z) = (n-1)!\, 2^{n-1}
\int_{\mathbb{H}^n} f(w,s)\,
\frac{J_{n-1}(\sqrt{\tau}|z-w|)}{(\sqrt{\tau}|z-w|)^{n-1}}\, dw\, ds,
\qquad \tau \ge 0.
\end{equation*}

Fix $\tau>0$ and let $k \to \infty$, $\lambda \to 0$ along $(2k+n)|\lambda|
= \tau$. Two degenerations occur simultaneously.

\begin{enumerate}
\item[(a)] The phase $e^{-i\lambda/2\operatorname{Im}(z\cdot\bar w)}$
defining $*_{-\lambda}$ tends to $1$, so twisted convolution degenerates
to ordinary convolution $(u*v)(z) = \int_{\mathbb{C}^n} u(z-w)v(w)\, dw$
on $\mathbb{C}^n \cong \mathbb{R}^{2n}$.

\item[(b)] By the Mehler-Heine-type asymptotics for Laguerre
polynomials \cite[Ch.~VIII]{Szego1975},
\begin{equation*}
\frac{\varphi^{n-1}_{k,\lambda}(z)}{\varphi^{n-1}_{k,\lambda}(0)} \longrightarrow
K_\tau(z) := (n-1)!\, 2^{n-1}\,
\frac{J_{n-1}(\sqrt{\tau}|z|)}{(\sqrt{\tau}|z|)^{n-1}},
\end{equation*}
uniformly on compact sets, as $k \to \infty$, $\lambda \to 0$,
$(2k+n)|\lambda| \to \tau$.
\end{enumerate}

Note $K_\tau(0) = 1$ for every $\tau > 0$, and
\begin{equation*}
K_\tau(z) = \int_{S^{2n-1}} e^{i\sqrt{\tau}\,\omega\cdot z}\, d\omega,
\end{equation*}
the Fourier transform of normalised surface measure on the sphere of
radius $\sqrt{\tau}$ in $\mathbb{R}^{2n}$.

Combining (a) and (b), $W_{S,\lambda}(g) = g *_{-\lambda}
\varphi^{n-1}_{k,\lambda}$ converges, up to this normalisation, to
ordinary convolution against the Bessel kernel $K_\tau$. We call
\begin{equation*}
W_{S,\infty,\tau}(g) := g * K_\tau
\end{equation*}
the \textit{Strichartz Weyl transform at the boundary}. Writing
$f^0(w) := \int_{\mathbb{R}} f(w,s)\, ds$,
\begin{align*}
\mathcal{F}_S(f)(0,\tau,z)
&= \int_{\mathbb{H}^n} f(w,s)\, K_\tau(z-w)\, dw\, ds \\
&= \big(f^0 * K_\tau\big)(z) \\
&= W_{S,\infty,\tau}\big(f^0\big)(z),
\end{align*}
the exact analogue of the relation between $\mathcal{F}_S$ and
$W_{S,\lambda}$ at interior points.

\begin{rmk}
By \Cref{prop:idempotent-self-adjoint},  $W_{S,\lambda}$ for $\lambda \ne 0$ is, up to the scalar
$(2\pi)^n|\lambda|^{-n}$, an orthogonal projection onto a
discrete eigenspace. In contrast, $W_{S,\infty,\tau}$ is not idempotent. It realises
the continuous spectral density $dE_\tau/d\tau$ of the ordinary
Laplacian $-\Delta$ on $\mathbb{C}^n \cong \mathbb{R}^{2n}$. The map $u
\mapsto u * K_\tau$ restricts $\widehat u$ to the sphere $|\xi| =
\sqrt{\tau}$ and reconstructs from it. Equivalently, $u * K_\tau$ is given by the Hankel transform of the
spherical mean of $u$, evaluated at $\sqrt{\tau}$ and re-expanded in
$z$.
\end{rmk}

\subsection{Riemann-Lebesgue lemmas and continuity on $\Omega$}
\label{sec:rl-continuity}

Throughout this subsection we use the normalised transform
$\widetilde W_{S,\lambda,k}$, together with the following facts,
recorded here for reference.

\begin{itemize}
\item \textit{Uniform bound:}
\[
|\widetilde W_{S,\lambda,k}(g)(z)| \le \|g\|_{L^1(\mathbb{C}^n)},
\qquad g \in L^1(\mathbb{C}^n),\ z \in \mathbb{C}^n.
\]

\item \textit{Mehler-Heine asymptotic and its error term}
\cite[Ch.~VIII]{Szego1975}:
\[
\widetilde\varphi^{n-1}_{k,\lambda}(z) = K_{(2k+n)|\lambda|}(z) +
m_k\big(\sqrt{|\lambda|}\,|z|\big), \qquad
|m_k(t)| \le C\, (2k+n)^{-(n-1)/2-3/4}, \quad 0 < t \le \omega.
\]
The error bound holds uniformly over the Laguerre argument in any
fixed range $0 < t \le \omega$, not merely as $t \to 0$.

\item \textit{Bessel bound} \cite[(1.71.6)]{Szego1975}:
\[
|J_{n-1}(t)| \le c_n\, t^{n-1}(1+t)^{-1/2}, \qquad t \ge 0.
\]

\item \textit{Definition and bound for $K_\tau$:}
\[
K_\tau(z) = \int_{S^{2n-1}} e^{i\sqrt{\tau}\,\omega\cdot z}\, d\omega,
\qquad |K_\tau(z)| \le 1.
\]
\end{itemize}

\subsubsection{Decay in $k$}

\begin{thm}\label{thm:decay-in-k}
Fix $\lambda \in \mathbb{R}^*$ and $z \in \mathbb{C}^n$. For every $g
\in L^1(\mathbb{C}^n)$,
\begin{equation*}
\widetilde W_{S,\lambda,k}(g)(z) \longrightarrow 0 \qquad \text{as } k \to \infty.
\end{equation*}
\end{thm}

\begin{proof}
By density and the uniform bound, it suffices to prove the claim for
$g \in C^\infty_c(\mathbb{C}^n)$, supported in $|w| \le R$. On this
compact set, the Mehler-Heine asymptotic gives, uniformly for $|w|
\le R$,
\begin{equation*}
\widetilde\varphi^{n-1}_{k,\lambda}(w) = K_{(2k+n)|\lambda|}(w) +
O\big((2k+n)^{-(n-1)/2-3/4}\big).
\end{equation*}

Fix $\delta > 0$.
\begin{align*}
\widetilde W_{S,\lambda,k}(g)(z)
&= \int_{\mathbb{C}^n} g(z-w)\, \widetilde\varphi^{n-1}_{k,\lambda}(w)\,
e^{-i\frac{\lambda}{2}\operatorname{Im}(z\cdot\bar w)}\, dw, \\[4pt]
|\widetilde W_{S,\lambda,k}(g)(z)|
&\le \int_{\mathbb{C}^n} |g(z-w)|\, \big|\widetilde\varphi^{n-1}_{k,\lambda}(w)\big|\, dw \\[4pt]
&= \int_{|w|<\delta} |g(z-w)|\, dw
\;+\; \int_{|w|\ge\delta} |g(z-w)|\,\big|\widetilde\varphi^{n-1}_{k,\lambda}(w)\big|\, dw.
\end{align*}

By absolute continuity of $\int|g|$, choosing $\delta$ small enough
makes the first term as small as we like, uniformly in $k$. For the
second term, the Mehler-Heine asymptotic gives
\begin{equation*}
\int_{|w|\ge\delta} |g(z-w)|\,\big|\widetilde\varphi^{n-1}_{k,\lambda}(w)\big|\, dw
\le \sup_{|w|\ge\delta} \big|K_{(2k+n)|\lambda|}(w)\big| \int_{\mathbb{C}^n} |g| \;+\;
O\big((2k+n)^{-(n-1)/2-3/4}\big).
\end{equation*}
As $(2k+n)|\lambda| \to \infty$, $K_{(2k+n)|\lambda|}(w) \to 0$
uniformly for $|w| \ge \delta$, so this term vanishes as $k \to
\infty$. The error term vanishes as $k \to \infty$ also. Letting $k
\to \infty$ and then $\delta \to 0$ gives $\widetilde
W_{S,\lambda,k}(g)(z) \to 0$.
\end{proof}

\begin{rmk}
By the fixed-$\lambda$ $L^1 \to L^\infty$
bound for Laguerre projections \cite{KochRicci2007}, together with
the duality $\|P_k\|_{1\to\infty} = \|P_k\|_{2\to\infty}^2$, there are
constants $c,C>0$, depending on $n$, such that
\begin{equation*}
c\,k^{n-1} \le \|W_{S,\lambda,k}\|_{1\to\infty} \le C\,k^{n-1}.
\end{equation*}
\end{rmk}

\subsubsection{Decay in $\lambda$, fixed $k$}

\begin{lem}
Fix $k \in \mathbb{N}$. There is a constant $C = C(n,k)$ such that
\begin{equation*}
\big|\widetilde\varphi^{n-1}_{k,\lambda}(w)\big| \le C\big(1+|\lambda||w|^2\big)^k
e^{-\frac{1}{4}|\lambda||w|^2}, \qquad w \in \mathbb{C}^n,\ \lambda \in \mathbb{R}^*.
\end{equation*}
Consequently, for every fixed $\delta > 0$,
\begin{equation*}
\sup_{|w| \ge \delta} \big|\widetilde\varphi^{n-1}_{k,\lambda}(w)\big| \longrightarrow 0
\qquad \text{as } |\lambda| \to \infty.
\end{equation*}
\end{lem}

\begin{proof}
Since $L^{n-1}_k(x)$ is a polynomial of fixed degree $k$ in $x \ge 0$,
there is a constant $C_1 = C_1(n,k)$ such that
\begin{equation*}
|L^{n-1}_k(x)| \le C_1(1+x)^k, \qquad x \ge 0.
\end{equation*}
Setting $x = \frac{1}{2}|\lambda||w|^2$,
\begin{equation*}
\big|\varphi^{n-1}_{k,\lambda}(w)\big| =
\Big|L^{n-1}_k\big(\tfrac12|\lambda||w|^2\big)\Big|\,
e^{-\frac14|\lambda||w|^2} \le
C_1\big(1+\tfrac12|\lambda||w|^2\big)^k e^{-\frac14|\lambda||w|^2},
\end{equation*}
and dividing by $L^{n-1}_k(0)$ gives the stated bound on
$\widetilde\varphi^{n-1}_{k,\lambda}$, with $C := C_1\cdot
2^k/L^{n-1}_k(0)$.

Taking $u := |\lambda||w|^2$, the bound reads
\begin{equation*}
\big|\widetilde\varphi^{n-1}_{k,\lambda}(w)\big| \le C(1+u)^k
e^{-u/4} =: h(u).
\end{equation*}
Since $h(u) \to 0$ as $u \to \infty$, and for $|w| \ge \delta$,
\begin{equation*}
u = |\lambda||w|^2 \ge |\lambda|\delta^2 \longrightarrow \infty
\qquad \text{as } |\lambda| \to \infty.
\end{equation*}
Thus we get $\sup_{|w|\ge\delta} \big|\widetilde\varphi^{n-1}_{k,\lambda}(w)\big|
\le h(|\lambda|\delta^2) \to 0$ as $|\lambda| \to \infty$.
\end{proof}

\begin{thm}
Fix $k \in \mathbb{N}$. For every $g \in L^1(\mathbb{C}^n)$,
\begin{equation*}
\|\widetilde W_{S,\lambda,k}(g)\|_{L^\infty(\mathbb{C}^n)} \longrightarrow 0
\qquad \text{as } |\lambda| \to \infty.
\end{equation*}
\end{thm}

\begin{proof}
Fix $\varepsilon>0$. By uniform absolute continuity, choose
$\delta>0$ with
\begin{equation*}
\sup_z \int_{B(z,\delta)} |g| < \varepsilon.
\end{equation*}
We have
\begin{align*}
|\widetilde W_{S,\lambda,k}(g)(z)|
&\le \int_{|w|<\delta} |g(z-w)|\, dw
\;+\; \int_{|w|\ge\delta} |g(z-w)|\,\big|\widetilde\varphi^{n-1}_{k,\lambda}(w)\big|\, dw \\
&\le \int_{|w|<\delta} |g(z-w)|\, dw
\;+\; \|g\|_1 \sup_{|\zeta|\ge\delta} \big|\widetilde\varphi^{n-1}_{k,\lambda}(\zeta)\big|.
\end{align*}
The first term is bounded by $\varepsilon$, uniformly in $z$. By the
preceding lemma, the second term tends to $0$ as $|\lambda| \to
\infty$, uniformly in $z$. Hence $\|\widetilde
W_{S,\lambda,k}(g)\|_{L^\infty(\mathbb{C}^n)} \to 0$ as $|\lambda| \to
\infty$.
\end{proof}

\subsubsection{Decay at the boundary ray}

\begin{thm} \label{thm:decay-boundary}
For every $g \in L^1(\mathbb{C}^n)$,
\begin{equation*}
\|g * K_\tau\|_{L^\infty(\mathbb{C}^n)} \longrightarrow 0 \qquad \text{as } \tau \to \infty.
\end{equation*}
\end{thm}

\begin{proof}
By density, we may take $g \in C_c(\mathbb{C}^n)$. Given
$\varepsilon>0$, absolute continuity of $\int g$ provides $\delta>0$
with
\begin{equation*}
\sup_z \int_{B(z,\delta)} |g| < \varepsilon.
\end{equation*}
We have
\begin{align*}
|(g*K_\tau)(z)|
&\le \int_{|w|<\delta} |g(z-w)|\, dw
\;+\; \int_{|w|\ge\delta} |g(z-w)|\,|K_\tau(w)|\, dw \\
&\le \varepsilon \;+\; \|g\|_1 \sup_{|\zeta|\ge\delta} |K_\tau(\zeta)|.
\end{align*}
As $\tau \to \infty$, the Bessel bound forces $K_\tau(\zeta) \to 0$
uniformly for $|\zeta| \ge \delta$, and the claim follows.
\end{proof}

\begin{rmk}
The Riemann-Lebesgue lemma for Hankel transforms is a classical fact, drawn on in [23]. The argument above supplies a self-contained proof of it, from the stated Bessel and kernel bounds.
\end{rmk}

\subsubsection{Riemann-Lebesgue lemma for $\mathcal{F}_S$ on $\Omega$}

\begin{thm}\label{thm:rlsft}
For every $f \in L^1(\mathbb{H}^n)$ and fixed $z \in \mathbb{C}^n$,
\begin{equation*}
\widetilde{\mathcal{F}_S}(f)(a,z) \longrightarrow 0 \qquad \text{as } |a| \to \infty \text{ in } \Omega.
\end{equation*}
\end{thm}

\begin{proof}
Let $a_j \to \infty$ in $\Omega$. Three cases arise.

\textbf{Case 1:} $a_j = (\lambda_j,k_j)$ with $|\lambda_j| \to \infty$.

By Tonelli's theorem, since $f \in L^1(\mathbb{H}^n)$, the function
$F(w) := \int_{\mathbb{R}} |f(w,s)|\, ds$ satisfies $F \in
L^1(\mathbb{C}^n)$, and in particular $s \mapsto f(w,s)$ is integrable
for a.e.\ $w$. For such $w$, the classical Riemann-Lebesgue lemma
gives
\begin{equation*}
f^{-\lambda_j}(w) \longrightarrow 0 \qquad \text{as } |\lambda_j| \to \infty.
\end{equation*}
Since $|f^{-\lambda_j}| \le F \in L^1(\mathbb{C}^n)$ for every $j$,
dominated convergence gives
\begin{equation*}
\|f^{-\lambda_j}\|_{L^1(\mathbb{C}^n)} \longrightarrow 0.
\end{equation*}
By the uniform bound,
\begin{equation*}
|\widetilde{\mathcal{F}_S}(f)(a_j,z)| \le
\|f^{-\lambda_j}\|_{L^1(\mathbb{C}^n)} \longrightarrow 0.
\end{equation*}

\textbf{Case 2:} $a_j = (\lambda_j,k_j)$, with $\lambda_j$ bounded and
$k_j \to \infty$.

Then necessarily $(2k_j+n)|\lambda_j| \to \infty$. Write $F(w) :=
\int_{\mathbb{R}} |f(w,s)|\, ds \in L^1(\mathbb{C}^n)$, so
$|f^{-\lambda_j}(w)| \le F(w)$ for every $j$. Fix $\varepsilon > 0$.
Since $F \in L^1(\mathbb{C}^n)$, choose $0 < \delta < R$ so that
\begin{equation*}
\int_{|w|<\delta} F(z-w)\, dw < \varepsilon, \qquad
\int_{|w|>R} F(z-w)\, dw < \varepsilon.
\end{equation*}
Splitting the domain of integration into $|w|<\delta$,
$\delta\le|w|\le R$, and $|w|>R$, and using
$|\widetilde\varphi^{n-1}_{k,\lambda}| \le 1$ on the outer two pieces,
\begin{align*}
|\widetilde{\mathcal{F}_S}(f)(a_j,z)|
&\le \int_{|w|<\delta} F(z-w)\, dw \\
&\quad + \sup_{\delta\le|w|\le R} \big|\widetilde\varphi^{n-1}_{k_j,\lambda_j}(w)\big|
\int_{\delta\le|w|\le R} F(z-w)\, dw \\
&\quad + \int_{|w|>R} F(z-w)\, dw.
\end{align*}
On $\delta \le |w| \le R$, since $|\lambda_j| \le M$ for some $M$, the
Mehler-Heine argument satisfies $\sqrt{|\lambda_j|}\,|w| \le
\sqrt{M}\,R =: \omega$, so the uniform-range error bound applies:
\begin{equation*}
\widetilde\varphi^{n-1}_{k_j,\lambda_j}(w) =
K_{(2k_j+n)|\lambda_j|}(w) + m_{k_j}\big(\sqrt{|\lambda_j|}\,|w|\big),
\qquad
\big|m_{k_j}\big(\sqrt{|\lambda_j|}\,|w|\big)\big| \le
C\,(2k_j+n)^{-(n-1)/2-3/4}.
\end{equation*}
As $k_j \to \infty$, the error term vanishes. By the Bessel bound,
$K_\tau(w) \to 0$ uniformly for $|w| \ge \delta$ as $\tau \to \infty$,
so $K_{(2k_j+n)|\lambda_j|}(w) \to 0$ uniformly on $\delta \le |w|
\le R$ as $j \to \infty$. Hence
\begin{equation*}
\sup_{\delta\le|w|\le R} \big|\widetilde\varphi^{n-1}_{k_j,\lambda_j}(w)\big| \longrightarrow 0,
\end{equation*}
and the middle term vanishes. Letting $j \to \infty$,
\begin{equation*}
\limsup_{j\to\infty} |\widetilde{\mathcal{F}_S}(f)(a_j,z)| \le
\varepsilon + 0 + \varepsilon = 2\varepsilon.
\end{equation*}
Since $\varepsilon > 0$ was arbitrary,
$\widetilde{\mathcal{F}_S}(f)(a_j,z) \to 0$.

\textbf{Case 3:} $a_j = (0,\tau_j) \in R_\infty$ with $\tau_j \to
\infty$.

Recall from \Cref{sec:boundary-ray} the identity $\widetilde{\mathcal{F}_S}(f)(0,\tau,z) =
(f^0*K_\tau)(z)$, where $f^0(w) := \int_{\mathbb{R}} f(w,s)\, ds$.
Thus
\begin{equation*}
\widetilde{\mathcal{F}_S}(f)(a_j,z) = (f^0*K_{\tau_j})(z).
\end{equation*}
Since $f^0 \in L^1(\mathbb{C}^n)$, \Cref{thm:decay-boundary}
gives
\begin{equation*}
\|f^0*K_\tau\|_{L^\infty(\mathbb{C}^n)} \to 0 \qquad \text{as }
\tau \to \infty.
\end{equation*}
Hence
\begin{equation*}
|\widetilde{\mathcal{F}_S}(f)(a_j,z)| \le
\|f^0*K_{\tau_j}\|_{L^\infty(\mathbb{C}^n)} \longrightarrow 0.
\end{equation*}
\end{proof}

\begin{rmk}
This is the same statement as \cite[Thm.~7.1]{Thangavelu2023}. The
proof given there establishes all three cases within a single
argument, using the Laguerre and Bessel asymptotics as intermediate
steps along the way. Here, these same asymptotic facts have already
been extracted and proved separately as standalone properties of
$\widetilde W_{S,\lambda,k}$ and $W_{S,\infty,\tau}$
(\Cref{thm:decay-in-k} and \Cref{thm:decay-boundary}), so the proof above
reduces to citing them case by case, rather than re-deriving the
estimates at the point of use.
\end{rmk}

\subsubsection{Continuity of $\mathcal{F}_S$ on $\Omega$}

Fix $z \in \mathbb{C}^n$ throughout this subsection.

\begin{lem}[Continuity in $\lambda$, fixed $k$]
For fixed $k$, the map $\lambda \mapsto \widetilde{\mathcal{F}_S}(f)(\lambda,k,z)$
is continuous on $\mathbb{R}^*$.
\end{lem}

\begin{proof}
Fix $\lambda_0 \in \mathbb{R}^*$. By $L^1$-continuity of translation
and modulation in the $t$-variable,
\begin{equation*}
\|f^{-\lambda} - f^{-\lambda_0}\|_{L^1(\mathbb{C}^n)} \longrightarrow 0
\qquad \text{as } \lambda \to \lambda_0.
\end{equation*}
Writing $\widetilde{\mathcal{F}_S}(f)(\lambda,k,z) =
\widetilde W_{S,\lambda,k}(f^{-\lambda})(z)$ and using the triangle
inequality,
\begin{align*}
\big|\widetilde{\mathcal{F}_S}(f)(\lambda,k,z) -
\widetilde{\mathcal{F}_S}(f)(\lambda_0,k,z)\big|
&\le \big|\widetilde W_{S,\lambda,k}(f^{-\lambda}-f^{-\lambda_0})(z)\big| \\
&\quad + \big|\widetilde W_{S,\lambda,k}(f^{-\lambda_0})(z) -
\widetilde W_{S,\lambda_0,k}(f^{-\lambda_0})(z)\big|.
\end{align*}
By the uniform bound, the first term is at most
$\|f^{-\lambda}-f^{-\lambda_0}\|_{L^1(\mathbb{C}^n)} \to 0$.\\
For the second term, the kernel $\widetilde\varphi^{n-1}_{k,\lambda}$ depends
smoothly on $\lambda$ at fixed $k$. So dominated convergence gives
\[
\widetilde W_{S,\lambda,k}(f^{-\lambda_0})(z) \to
\widetilde W_{S,\lambda_0,k}(f^{-\lambda_0})(z)
\qquad \text{as } \lambda \to \lambda_0.
\]
Hence $\widetilde{\mathcal{F}_S}(f)(\lambda,k,z) \to
\widetilde{\mathcal{F}_S}(f)(\lambda_0,k,z)$.
\end{proof}

\begin{lem}[Continuity in $\tau$ on $R_\infty$]
Fix $z$. The map $\tau \mapsto (g*K_\tau)(z)$ is uniformly continuous
on $(0,\infty)$ for $g \in L^1(\mathbb{C}^n)$.
\end{lem}

\begin{proof}
Write $\Psi(\tau) := \int_{S^{2n-1}} \widehat{g_z}(-\sqrt{\tau}\,\omega)\, d\omega$,
where $g_z(w) := g(z-w)$, so that $(g*K_\tau)(z) = \Psi(\tau)$. Since
$\widehat{g_z}$ is uniformly continuous on $\mathbb{C}^n$
\cite[Ch.~I]{SteinWeissIntro1971}, and
\begin{equation*}
\big|\sqrt{\tau}\,\omega - \sqrt{\tau_0}\,\omega\big| = \big|\sqrt{\tau} - \sqrt{\tau_0}\big|
\end{equation*}
uniformly in $\omega \in S^{2n-1}$, it follows that
$|\Psi(\tau) - \Psi(\tau_0)| \to 0$ as $\tau \to \tau_0$, uniformly in $z$.
\end{proof}

\begin{lem}[Boundary limit]
Let $a_j = (\lambda_j,k_j) \in R_{k_j}$ with $a_j \to (0,\tau)$, $\tau >
0$. Then $\widetilde{\mathcal{F}_S}(f)(a_j,z) \to
\widetilde{\mathcal{F}_S}(f)((0,\tau),z)$.
\end{lem}

\begin{proof}
Write $g_j := \widetilde\varphi^{n-1}_{k_j,\lambda_j}$, so that
\begin{equation*}
\widetilde{\mathcal{F}_S}(f)(a_j,z) = (f^{-\lambda_j} * g_j)(z),
\qquad
\widetilde{\mathcal{F}_S}(f)((0,\tau),z) = (f^0 * K_\tau)(z).
\end{equation*}
By the triangle inequality,
\begin{equation*}
\big|\widetilde{\mathcal{F}_S}(f)(a_j,z) - \widetilde{\mathcal{F}_S}(f)((0,\tau),z)\big|
\le I + II,
\end{equation*}
where
\begin{equation*}
I := \big|(f^{-\lambda_j} * g_j)(z) - (f^0 * g_j)(z)\big|,
\qquad
II := \big|(f^0 * g_j)(z) - (f^0 * K_\tau)(z)\big|.
\end{equation*}

\textbf{Bound on $I$.} By the uniform bound, $|g_j(w)| \le 1$ for all $j$ and $w$, so
\begin{align*}
I &= \left|\int_{\mathbb{C}^n} \big(f^{-\lambda_j}(z-w) - f^0(z-w)\big) g_j(w)\, dw\right| \\
&\le \int_{\mathbb{C}^n} |f^{-\lambda_j}(z-w) - f^0(z-w)|\, dw \\
&= \|f^{-\lambda_j} - f^0\|_{L^1} \to 0.
\end{align*}

\textbf{Bound on $II$.} Here
\begin{equation*}
II = \left|\int_{\mathbb{C}^n} f^0(z-w)\big(g_j(w) - K_\tau(w)\big)\, dw\right|.
\end{equation*}
Pointwise, $g_j(w) \to K_\tau(w)$ as $j \to \infty$ (Mehler-Heine convergence, using the
error-term bound above). By the uniform bound and $|K_\tau(z)| \le 1$, the integrand is dominated as
\begin{equation*}
\big|f^0(z-w)\big(g_j(w) - K_\tau(w)\big)\big| \le 2\,|f^0(z-w)|,
\end{equation*}
and $w \mapsto f^0(z-w) \in L^1(\mathbb{C}^n)$. Dominated convergence theorem gives $II \to 0$.

Combining the two bounds, $I + II \to 0$, so
$\widetilde{\mathcal{F}_S}(f)(a_j,z) \to \widetilde{\mathcal{F}_S}(f)((0,\tau),z)$.
\end{proof}

\begin{thm}
For fixed $z$, the map $a \mapsto \widetilde{\mathcal{F}_S}(f)(a,z)$ is
continuous on $\Omega$.
\end{thm}

\begin{proof}
Let $a_0 \in \Omega$ and let $(a_j)$ be any sequence in $\Omega$ with $a_j \to a_0$.

\textbf{Case 1: $a_0 \in R_k$ for some $k$.} Since $R_k$ is open in $\Omega$,
$a_j \in R_k$ for all large $j$, and $a_j \to a_0$ means $\lambda_j \to
\lambda_0$ with $k$ fixed. Continuity now follows directly from the
continuity-in-$\lambda$ lemma.

\textbf{Case 2: $a_0 = (0,\tau_0) \in R_\infty$.} Split $(a_j)$ into two
subsequences according to whether $a_j \in R_\infty$ or $a_j \in R_{k_j}$
for some finite $k_j$.

\begin{itemize}
\item For the subsequence with $a_j \in R_\infty$, write $a_j = (0,\tau_j)$
with $\tau_j \to \tau_0$. Continuity along this subsequence is exactly
the continuity-in-$\tau$ lemma.

\item For the subsequence with $a_j = (\lambda_j,k_j) \in R_{k_j}$,
$a_j \to (0,\tau_0)$ forces $\lambda_j \to 0$. Continuity along this
subsequence is same as the boundary-limit lemma.
\end{itemize}

Since every subsequence of $(a_j)$ has a further subsequence along which
$\widetilde{\mathcal{F}_S}(f)(a_j,z) \to \widetilde{\mathcal{F}_S}(f)(a_0,z)$
by one of the two cases above, and the two cases are exhaustive, the full
sequence converges: $\widetilde{\mathcal{F}_S}(f)(a_j,z) \to
\widetilde{\mathcal{F}_S}(f)(a_0,z)$.
\end{proof}

\begin{cor}
For fixed $z$, $a \mapsto \widetilde{\mathcal{F}_S}(f)(a,z)$ is
uniformly continuous on $\Omega$.
\end{cor}

\begin{proof}
Let $\varepsilon > 0$. By the Riemann-Lebesgue lemma (\Cref{thm:rlsft}), $\widetilde{\mathcal{F}_S}(f)(a,z) \to 0$ as $|a| \to \infty$
in $\Omega$. So there is $N > 0$ such that
\begin{equation*}
|a| \ge N \implies \big|\widetilde{\mathcal{F}_S}(f)(a,z)\big| < \varepsilon/2.
\end{equation*}
In particular, for any $a, a' \in \Omega$ with $|a|, |a'| \ge N$,
\begin{equation*}
\big|\widetilde{\mathcal{F}_S}(f)(a,z) - \widetilde{\mathcal{F}_S}(f)(a',z)\big|
\le \big|\widetilde{\mathcal{F}_S}(f)(a,z)\big| + \big|\widetilde{\mathcal{F}_S}(f)(a',z)\big|
< \varepsilon.
\end{equation*}

On the compact set $K_N := \{a \in \Omega : |a| \le N\}$, the preceding
theorem gives continuity of $a \mapsto \widetilde{\mathcal{F}_S}(f)(a,z)$
at every point. In particular, $\widetilde{\mathcal{F}_S}(f)(a,z)$ is uniformly continuous on $K_N$. So there is $\delta > 0$ such that
\begin{equation*}
a, a' \in K_N, \ |a-a'| < \delta \implies
\big|\widetilde{\mathcal{F}_S}(f)(a,z) - \widetilde{\mathcal{F}_S}(f)(a',z)\big| < \varepsilon.
\end{equation*}

Set $\delta' := \min(\delta, 1)$, and replace $N$ by $N+1$ throughout the compactness
step above, so that the uniform continuity bound holds on $K_{N+1}$ instead of $K_N$. The decay bound
established for $|a| \ge N$ remains valid for $|a| \ge N+1$, so this
substitution is free.

Let $a, a' \in \Omega$ with $|a - a'| < \delta'$. Two cases arise.

\textbf{Case 1: $|a|, |a'| \ge N$.} The decay bound applies to both $a$
and $a'$ directly, giving
\begin{equation*}
\big|\widetilde{\mathcal{F}_S}(f)(a,z) - \widetilde{\mathcal{F}_S}(f)(a',z)\big|
\le \big|\widetilde{\mathcal{F}_S}(f)(a,z)\big| + \big|\widetilde{\mathcal{F}_S}(f)(a',z)\big|
< \varepsilon.
\end{equation*}

\textbf{Case 2: at least one of $|a|, |a'|$ is $< N$.} Say $|a| < N$.
Since $\delta' \le 1$,
\begin{equation*}
|a'| \le |a| + |a-a'| < N + 1.
\end{equation*}
So both $a, a' \in K_{N+1}$, and the uniform continuity bound on
$K_{N+1}$ gives
\begin{equation*}
\big|\widetilde{\mathcal{F}_S}(f)(a,z) - \widetilde{\mathcal{F}_S}(f)(a',z)\big| < \varepsilon.
\end{equation*}

In either case, $a, a' \in \Omega$ with $|a-a'| < \delta'$ satisfy
\begin{equation*}
\big|\widetilde{\mathcal{F}_S}(f)(a,z) - \widetilde{\mathcal{F}_S}(f)(a',z)\big| < \varepsilon,
\end{equation*}
proving uniform continuity.
\end{proof}

\begin{lem}[Continuity in $z$, fixed $a$]
For fixed $a = (\lambda,k)$, the map $z \mapsto
\widetilde{\mathcal{F}_S}(f)(a,z)$ is continuous on $\mathbb{C}^n$.
\end{lem}

\begin{proof}
Write $\varphi := \widetilde\varphi^{n-1}_{k,\lambda}$, so that
$\widetilde{\mathcal{F}_S}(f)(a,z) = (f^{-\lambda} * \varphi)(z)$. Fix
$z_0 \in \mathbb{C}^n$ and $\varepsilon > 0$.

Since $f^{-\lambda} \in L^1(\mathbb{C}^n)$, choose a continuous,
compactly supported $g$ with
\begin{equation*}
\|f^{-\lambda} - g\|_{L^1(\mathbb{C}^n)} < \varepsilon.
\end{equation*}
Split
\begin{equation*}
\big(f^{-\lambda}*\varphi\big)(z) - \big(f^{-\lambda}*\varphi\big)(z_0)
= I + II + III,
\end{equation*}
where
\begin{align*}
I &:= \int_{\mathbb{C}^n} \big(f^{-\lambda}(z-w) - g(z-w)\big)\varphi(w)\,dw, \\
II &:= \int_{\mathbb{C}^n} \big(g(z-w) - g(z_0-w)\big)\varphi(w)\,dw, \\
III &:= \int_{\mathbb{C}^n} \big(g(z_0-w) - f^{-\lambda}(z_0-w)\big)\varphi(w)\,dw.
\end{align*}

\textbf{Bounds on $I$ and $III$.} By the uniform bound $|\varphi(w)| \le 1$,
\begin{equation*}
|I| \le \int_{\mathbb{C}^n} |f^{-\lambda}(z-w) - g(z-w)|\,dw
= \|f^{-\lambda} - g\|_{L^1} < \varepsilon,
\end{equation*}
and likewise $|III| < \varepsilon$, using the translation-invariance of
Lebesgue measure.

\textbf{Bound on $II$.} Since $g$ is continuous with compact support,
say $\operatorname{supp} g \subset \{|w|\le R\}$, it is uniformly
continuous on $\mathbb{C}^n$. So there is $\eta > 0$ such that
\begin{equation*}
|z - z_0| < \eta \implies \sup_{w \in \mathbb{C}^n} |g(z-w) - g(z_0-w)| < \varepsilon.
\end{equation*}
For such $z$,
\begin{equation*}
|II| \le \int_{\mathbb{C}^n} |g(z-w) - g(z_0-w)| \, |\varphi(w)| \,dw
\le \varepsilon \int_{|w| \le R+|z-z_0|} 1 \,dw
\le \varepsilon \cdot C_R,
\end{equation*}
where $C_R$ is the volume of a fixed ball containing
$\operatorname{supp} g(z-\cdot)$ for all $z$ with $|z-z_0|<\eta$, and we
again used $|\varphi(w)|\le 1$.

Combining, for $|z-z_0| < \eta$,
\begin{equation*}
\big|\widetilde{\mathcal{F}_S}(f)(a,z) - \widetilde{\mathcal{F}_S}(f)(a,z_0)\big|
= |I+II+III| < \varepsilon(2 + C_R).
\end{equation*}
Since $\varepsilon > 0$ was arbitrary, $z \mapsto
\widetilde{\mathcal{F}_S}(f)(a,z)$ is continuous at $z_0$, and $z_0$ was
arbitrary.
\end{proof}

\begin{rmk}
Together, the preceding results of this subsection recover, in full, the
continuity and vanishing-at-infinity statements for $\mathcal{F}_S$
given in \cite[Thm.~7.1]{Thangavelu2023}. Both the original proof and
the one given above proceed by the same underlying three-case split
across $\Omega$. Here each case is instead handled through
$L^1$-continuity of $f^{-\lambda}$ or $f^0$, dominated convergence,
and the Mehler-Heine convergence of the Laguerre kernel to the Bessel
kernel $K_\tau$, assembled into separate lemmas rather than treated in
a single argument.
\end{rmk}

\subsubsection{Bochner-Riesz means on the Heisenberg fan}

\begin{defn} \label{def:br}
For $\eta \ge 0$, $R>0$, and $f \in L^2(\mathbb{H}^n)$, define
\begin{equation*}
S^\eta_R f(z,t) := \int_\Omega \Big(1-\frac{\tau(a)}{R^2}\Big)_+^\eta
\widehat f(a,w)\, e_a\big((-w,0)(z,t)\big)\, dw\, d\nu(a),
\end{equation*}
where
\begin{equation*}
\tau(a) := (2k+n)|\lambda| \text{ for } a \in R_k.
\end{equation*}
\end{defn}

\begin{prop}\label{prop:BR-operator-identity}
For $\eta \geq 0$, $R > 0$, and $f \in L^2(\mathbb{H}^n)$,
\[
S^\eta_R f = (1 - R^{-2}\mathcal{L})^\eta_+ f.
\]
\end{prop}

\begin{proof}
Fix $a = (\lambda, k, \varepsilon) \in \Omega$, and write $f_a := f * \widetilde\varphi^{\,n-1}_{k,\lambda}$
for the component of $f$ at the point $a$.
By the joint spectral theory of $\mathcal{L}$ and $T$ established by Strichartz \cite{Strichartz1991},
every such component satisfies
\begin{align*}
-\mathcal{L} f_a = \tau(a)\, f_a,
\end{align*}
where $\tau(a) = (2k+n)|\lambda|$, for every $f \in L^2(\mathbb{H}^n)$,
not merely when $f$ itself is a joint eigenfunction.

Consequently, for any polynomial $p$, the component of $p(\mathcal{L}) f$ at $a$ is $p(\tau(a)) f_a$,
and by density and a standard limiting argument this extends to any bounded Borel function $F$
on $[0, \infty)$,
\begin{align} \label{eq:bddborel}
(F(\mathcal{L}) f)_a = F(\tau(a))\, f_a \qquad \text{for every } a \in \Omega.
\end{align}

Now take $F(\tau) = (1 - R^{-2}\tau)^\eta_+$.
By \Cref{def:br}, the operator $S^\eta_R$ acts on the component $f_a$ by the scalar multiplier
$(1 - R^{-2}\tau(a))^\eta_+$,
\begin{align*}
(S^\eta_R f)_a = (1 - R^{-2}\tau(a))^\eta_+\, f_a.
\end{align*}

Combining the above equation with \Cref{eq:bddborel}, $(S^\eta_R f)_a = (F(\mathcal{L}) f)_a$ for every $a \in \Omega$,
and since the fan decomposition recovers $f$ from its components (\Cref{thm:inversion-omega}), this gives
\begin{align*}
S^\eta_R f = (1 - R^{-2}\mathcal{L})^\eta_+ f. 
\end{align*}
\end{proof}

\begin{thm}[Bochner-Riesz means for $\mathcal{F}_S$]
\label{thm:bochner-riesz}
Let $d = 2n+1$. In the sharp range $1 \le p \le \frac{2(d+1)}{d+4}$ and its
dual, $S^\eta_R$ is bounded on $L^p(\mathbb{H}^n)$, uniformly in $R$, for
$\eta > d|1/p-1/2|-1/2$. More generally, for $1 \le p < p_n$, where
$p_n \to 2$ as $n \to \infty$, the same uniform bound holds under the
weaker Sobolev-smoothness condition $\eta > d|1/p-1/2|$. In either range,
$S^\eta_R f \to f$ in $L^p(\mathbb{H}^n)$ as $R \to \infty$.
\end{thm}

\begin{proof}
By \Cref{prop:BR-operator-identity}, $S^\eta_R = (1-R^{-2}\mathcal{L})^\eta_+$,
so uniform boundedness of $S^\eta_R$ on $L^p(\mathbb{H}^n)$ is the same as the
classical Bochner-Riesz problem for the sub-Laplacian $\mathcal{L}$, whether
the Riesz means of order $\eta$ for $\mathcal{L}$ are bounded uniformly in
$R$.

We get the sharp bound in the range $1 \le p \le \frac{2(d+1)}{d+4}$ and its
dual, together with the wider range $1 \le p < p_n$ under the weaker
Sobolev-smoothness condition, from \cite{MullerNiedorfSeegerStovall2026}. Sharpness of the exponent
$\eta > d|1/p-1/2|-1/2$ in the first range is due to \cite{Martini2025}.
Mean convergence $S^\eta_R f \to f$ in $L^p(\mathbb{H}^n)$ then follows
from the uniform bound, by density.
\end{proof}

\begin{rmk}
This result is an immediate transcription once $\tau(a)$ is recognised as
the eigenvalue of $\mathcal{L}$. The statement above and the classical
Bochner-Riesz problem for $\mathcal{L}$ are the same statement in
different notation. The genuinely open part of this problem, the $L^p$
range beyond $p \le \frac{2(d+1)}{d+4}$ and its dual, remains open
precisely because it is open at the operator level, with nothing specific
to $\mathcal{F}_S$ contributing to or obstructing it.
\end{rmk}

\section{The Hecke-Bochner Formula and the Image of the Schwartz Space}
\label{sec:hecke-bochner-schwartz}

The material in this section has a different character from
\Cref{sec:new-results,sec:core-toolkit}. The Hecke-Bochner formula
itself (\Cref{sec:hecke-bochner}) is reached through
$W_{S,\lambda}$. The full image of $\mathcal{S}(\mathbb{H}^n)$, however, is not. It requires a
compatibility condition linking values on different rays as
$\lambda \to 0$, which is invisible to $W_{S,\lambda}$ at any single
fixed ray, and is instead reached by citing the harder, global
direction of a theorem of Astengo, Di Blasio, and Ricci
\cite{AstengoDiBlasioRicci2007}.

\subsection{Radial functions and the Gelfand transform}
\label{sec:gelfand}

A function $f$ on $\mathbb{H}^n$ is radial if $f(z,t) = f_0(|z|,t)$ for
a unique function $f_0$ on $\mathbb{R}_+\times\mathbb{R}$. For radial
$f$, the special Hermite expansion of $f^\lambda$ reduces to a
Laguerre expansion,
\begin{equation*}
f^\lambda(z) = (2\pi)^{-n}|\lambda|^n \sum_{k=0}^\infty
R^{n-1}_k(\lambda,f)\, \varphi^{n-1}_{k,\lambda}(z),
\end{equation*}
where $R^{n-1}_k(\lambda,f)$ is the $k$-th Laguerre coefficient,
\begin{equation*}
R^{n-1}_k(\lambda,f) = \frac{k!(n-1)!}{(k+n-1)!} \int_{\mathbb{C}^n}
f^\lambda(z)\, \varphi^{n-1}_{k,\lambda}(z)\, dz.
\end{equation*}
As elements of $L^1(\mathbb{H}^n/U(n))$ are precisely the $U(n)$-biinvariant
functions on $G_n$, this expansion shows that $(G_n,U(n))$ is a Gelfand
pair, and that $L^1(\mathbb{H}^n/U(n))$ is a commutative convolution
algebra. The multiplicative linear functionals on this algebra are given
by the bounded spherical functions
\begin{equation*}
\widetilde e_a(z,t) = \frac{k!(n-1)!}{(k+n-1)!}\, e^{n-1}_{k,\lambda}(z,t),
\qquad \chi_\tau(z,t) = (n-1)!\, 2^{n-1}\,
\frac{J_{n-1}(\sqrt{\tau}|z|)}{(\sqrt{\tau}|z|)^{n-1}}, \quad \tau \ge 0,
\end{equation*}
established in \cite{BensonJenkinsRatcliff1998}. The Gelfand spectrum of
$L^1(\mathbb{H}^n/U(n))$ is therefore the Heisenberg fan
$\Omega$. The Gelfand transform $\mathcal{G}f$ is defined by
\begin{equation*}
\mathcal{G}f(a) = \int_{\mathbb{H}^n} f(z,t)\, \widetilde e_a(z,t)\,
dz\, dt, \qquad \mathcal{G}f(0,\tau) = \int_{\mathbb{H}^n} f(z,t)\,
\chi_\tau(z,t)\, dz\, dt,
\end{equation*}
and, for radial $f$,
\begin{equation*}
\widehat f(a,z) = \mathcal{G}f(a)\, e^{n-1}_{k,\lambda}(z,0).
\end{equation*}

\subsection{The Hecke-Bochner formula}
\label{sec:hecke-bochner}

The pair $(G_n,U(n))$ is a Gelfand pair \cite{AstengoDiBlasioRicci2007}.
This means $L^1_{U(n)}(\mathbb{H}^n)$ is commutative under convolution, as
recorded in the preceding subsection. Given a pair $(p,q)$ of
non-negative integers, a polynomial $P(z,\bar z)$ on $\mathbb{C}^n$ is
a bigraded solid harmonic if it is harmonic and homogeneous of degree
$p$ in $z$ and $q$ in $\bar z$.

Given $p,q \ge 0$, let $W^{(n+p+q)}_{S,\lambda,k-q}$ and $W^{(n+p+q)}_{S,\lambda,k-p}$ denote the
analogous transform built from Laguerre functions of the shifted type $n+p+q$ in place of
$n-1$.\\
We set $W^{(n+p+q)}_{S,\lambda,k-q} := 0$ when $k<q$, and
$W^{(n+p+q)}_{S,\lambda,k-p} := 0$ when $k<p$.

\begin{thm}[Hecke-Bochner formula]
Let $(p,q)$ be a pair of non-negative integers, and let $g$ be a radial function on
$\mathbb{C}^n$ such that $gP \in L^1 \cap L^2(\mathbb{C}^n)$ for every solid harmonic $P$ of
bidegree $(p,q)$. Then, at the $W_{S,\lambda,k}$ level,
\begin{align*}
\lambda > 0:\quad & W_{S,\lambda,k}(gP) = (2\pi)^{-(p+q)}\lambda^{p+q}\, P
\cdot W^{(n+p+q)}_{S,\lambda,k-q}(g), \\
\lambda < 0:\quad & W_{S,\lambda,k}(gP) = (2\pi)^{-(p+q)}|\lambda|^{p+q}\, P
\cdot W^{(n+p+q)}_{S,\lambda,k-p}(g),
\end{align*}
\end{thm}

\begin{proof} We work on the two cases separately.\\
\noindent \textbf{Case 1:} $\lambda>0$. \\
This case is the $W_{S,\lambda}$-level restatement of the
Hecke-Bochner-type identity for $\mathcal{F}_S$ established in \cite[Thm.~7.6]{Thangavelu2023},
via the relation between $\mathcal{F}_S$ and $W_{S,\lambda}$ (Remark~\ref{rmk:comparing-FS}).
That identity is itself built on the identity of Geller \cite{Geller1980}.\\
Thus for $\lambda_0>0$,
\begin{equation}
  W_{S,\lambda_0,k}(gP) = (2\pi)^{-(p+q)}\lambda_0^{p+q}\, P \cdot
  W^{(n+p+q)}_{S,\lambda_0,k-q}(g). \label{eq:hb-case1}
\end{equation}

\noindent \textbf{Case 2:} $\lambda<0$.\\
Apply the conjugation proposition
(Proposition~\ref{prop:conjugation}) to $h = gP$ at $\mu=\lambda_0$. Since $g$ is real-valued,
$\overline{gP} = g\bar P$, and $\bar P$ has bidegree $(q,p)$. This, together with
\eqref{eq:hb-case1}, gives
\begin{equation}
  W_{S,-\lambda_0,k}(g\bar P) = \overline{W_{S,\lambda_0,k}(gP)}
  = (2\pi)^{-(p+q)}\lambda_0^{p+q}\, \bar P \cdot
  \overline{W^{(n+p+q)}_{S,\lambda_0,k-q}(g)}. \label{eq:hb-case2-mid}
\end{equation}

We apply the same conjugation proposition a second time, now to the inner factor alone at
the shifted dimension $n+p+q$. Using again that $g$ is real, so $\bar g = g$, this gives
\[
  \overline{W^{(n+p+q)}_{S,\lambda_0,k-q}(g)}
  = \overline{W^{(n+p+q)}_{S,\lambda_0,k-q}(\bar g)}
  = W^{(n+p+q)}_{S,-\lambda_0,k-q}(g).
\]
Substituting this into \eqref{eq:hb-case2-mid},
\[
  W_{S,-\lambda_0,k}(g\bar P) = (2\pi)^{-(p+q)}\lambda_0^{p+q}\, \bar P \cdot
  W^{(n+p+q)}_{S,-\lambda_0,k-q}(g).
\]
Writing $\lambda := -\lambda_0 < 0$ and $Q := \bar P$, a harmonic of bidegree $(q,p)$, this reads
\[
  W_{S,\lambda,k}(gQ) = (2\pi)^{-(p+q)}|\lambda|^{p+q}\, Q \cdot
  W^{(n+p+q)}_{S,\lambda,k-q}(g).
\]
Relabelling $Q$ as $P$, now understood to have bidegree $(p,q)$ in its own right, gives
\[
  W_{S,\lambda,k}(gP) = (2\pi)^{-(p+q)}|\lambda|^{p+q}\, P \cdot
  W^{(n+p+q)}_{S,\lambda,k-p}(g), \qquad \lambda<0,
\]
which is the stated identity, since conjugation is a bijection on bigraded solid harmonics
of each bidegree.
\end{proof}

The treatment of the case $\lambda<0$ given in \cite[Thm.~7.6]{Thangavelu2023} proceeds by
interchanging the roles of $p$ and $q$ directly. The derivation above instead obtains the
$\lambda<0$ case from the $\lambda>0$ case through the conjugation identity of
Proposition~\ref{prop:conjugation}, giving an independent route to the same statement.

\begin{rmk}
Transplantation theorems, in the classical sense, show that
$L^p$-boundedness properties of coefficient operators transfer, up to
constants, across different values of the expansion parameter, and so,
in effect, across different ambient dimensions
\cite{AskeyWainger1966}. This resembles the dimension shift $n \to
n+p+q$ appearing in the Hecke-Bochner formula just proved, but the
resemblance is only superficial. The Hecke-Bochner formula is an
exact algebraic identity, determined by the bidegree of a solid
harmonic, while classical transplantation theorems are genuine
$L^p$-norm inequalities across a continuously varying parameter,
proved by independent analytic estimates. Neither the Hecke-Bochner
formula nor classical transplantation theorems derive one another,
and we do not pursue the connection further here.
\end{rmk}

\subsection{The Schwartz-space isomorphism at a fixed ray}

Let $\mathcal{H}_k \subset L^2(\mathbb{C}^n)$ denote the range of $P_k$
and $\mathcal{S}_k(\mathbb{C}^n) := \mathcal{S}(\mathbb{C}^n) \cap
\mathcal{H}_k$.

\begin{thm}
\label{thm:schwartz-iso}
$W_{S,\lambda}$ maps $\mathcal{S}(\mathbb{C}^n)$ into $\mathcal{S}_k(\mathbb{C}^n)$. In
particular, $W_{S,\lambda}$ is a linear isomorphism of $\mathcal{S}_k(\mathbb{C}^n)$ onto
itself.
\end{thm}

\begin{proof}
By the self-reproducing range condition (Theorem~\ref{thm:self-reproducing}),
$W_{S,\lambda}$ maps $\mathcal{S}(\mathbb{C}^n)$ into $\mathcal{S}_k(\mathbb{C}^n)$.

For $h \in \mathcal{S}_k(\mathbb{C}^n)$, $P_kh = h$, so
\[
  W_{S,\lambda}(h) = (2\pi)^n|\lambda|^{-n}\, h,
\]
giving the isomorphism claim.

Take $g := (2\pi)^{-n}|\lambda|^n h$. Then $W_{S,\lambda}(g) = h$.

So every $h \in \mathcal{S}_k(\mathbb{C}^n)$ arises this way. This proves
$W_{S,\lambda}(\mathcal{S}(\mathbb{C}^n)) = \mathcal{S}_k(\mathbb{C}^n)$.
\end{proof}

\begin{rmk}
This is the fixed-$k$ piece of the Schwartz-space image theorem for $\mathcal{F}_S$ stated in \cite[Thms.~7.7-7.8]{Thangavelu2023}. This full theorem requires both a decay-in-$k$ characterisation within a ray and a compatibility condition linking values on different rays as $\lambda \to 0$. The latter cannot be formulated at the level of a single fixed ray.
\end{rmk}

\section{Limits of the Fixed-Ray Route}
\label{sec:limits}

Every result in \Cref{sec:core-toolkit,sec:new-results} was reached by
assembling the fixed-ray theory of $W_{S,\lambda}$ developed in
\Cref{sec:swt}. This was done by integrating over $\lambda$,
summing compatibly over $k$, or transplanting estimates ray by ray.
Two existing results concerning $\mathcal{F}_S$ do not admit this
route. We give both the same level of attention as the positive
results above, since a precise account of where a method stops is as
much a part of understanding it as an account of where it succeeds.

\subsection{Ingham's theorem for the Strichartz Fourier transform}
\label{sec:ingham}

\begin{thm}[Ingham's theorem for the group Fourier transform, \cite{BagchiGangulySarkarThangavelu2022}]
Let $\Phi$ be a nonnegative even function on $\mathbb{R}$ decreasing to
$0$ at infinity. There exists a compactly supported function $f$ on
$\mathbb{H}^n$ whose group Fourier transform satisfies the decay
condition
\begin{equation*}
|\widehat f(\lambda)| \le c\, e^{-\sqrt{H(\lambda)}}\, \Phi(\sqrt{H(\lambda)})
\end{equation*}
if and only if $\int_1^\infty \Phi(t)\, t^{-1}\, dt < \infty$, where
$H(\lambda)$ denotes the Hermite operator at parameter $\lambda$.
\end{thm}

This operator-valued statement, due to Bagchi, Ganguly, Sarkar, and
Thangavelu \cite{BagchiGangulySarkarThangavelu2022}, restates the
classical Euclidean Ingham theorem for the Hermite operator $H(\lambda)$
at each fixed $\lambda$. In \cite[Thm.~7.10]{Thangavelu2023}, Thangavelu reformulates it at the $\mathcal{F}_S$ level, directly from the group-level statement rather than from any fixed-ray fact about $W_{S,\lambda}$, and we record that reformulation here.

\begin{thm}[Ingham's theorem for $\mathcal{F}_S$] (\cite{Thangavelu2023})
\label{thm:ingham}
Let $\Phi$ be a nonnegative even function on $\mathbb{R}$ decreasing to
$0$ at infinity. There exists a compactly supported function $f$ on
$\mathbb{H}^n$ whose Strichartz Fourier transform satisfies the decay
condition
\begin{equation*}
|\lambda|^n \int_{\mathbb{C}^n} |\widehat f(a,z)|^2\, dz \le c\,
\frac{(k+n-1)!}{k!(n-1)!}\, e^{-\sqrt{(2k+n)|\lambda|}}\,
\Phi\big(\sqrt{(2k+n)|\lambda|}\big)
\end{equation*}
if and only if $\int_1^\infty \Phi(t)\, t^{-1}\, dt < \infty$.
\end{thm}

\begin{proof}(\cite{Thangavelu2023})
This is obtained from the operator-valued form of Ingham's theorem via
the identity
\begin{equation*}
\int_{\mathbb{C}^n} |\widehat f(a,z)|^2\, dz = (2\pi)^n|\lambda|^{-n}\,
\|\widehat f(-\lambda)A_k^{(\lambda)}\|_{HS}^2,
\end{equation*}
obtained by combining the relation between the Weyl transform and
twisted convolution with the Plancherel identity for the Weyl
transform. The proof of the operator-valued theorem given in
\cite{BagchiGangulySarkarThangavelu2022} is modified by substituting
this identity for the Hilbert-Schmidt norm of $\widehat
f(-\lambda)A_k^{(\lambda)}$. Since $H(\mu) = -\Delta + \mu^2|x|^2$ depends on $\mu$ only through
$\mu^2$ and $A_k^{(\lambda)}$ depends on $\lambda$ only through
$|\lambda|$ (\Cref{sec:motiongroup}), the hypothesis and conclusion
of the operator-valued theorem coincide at $\mu=-\lambda$ and
$\mu=\lambda$, so this substitution is valid.
\end{proof}

\begin{rmk}[Why no fixed-ray route exists]
Theorem~\ref{thm:ingham} is an ``if and only if,'' and its two
directions place different demands on the fixed-ray method. The
necessity direction, that divergence of $\int_1^\infty
\Phi(t)t^{-1}\,dt$ forces any witness to vanish, reduces, in
\cite{BagchiGangulySarkarThangavelu2022}, to a fixed-$\lambda$
statement: an analogue of Chernoff's theorem for the single special
Hermite operator $L_\lambda$, applied at each $\lambda$ separately and
then assembled over $\lambda$. This is the same shape as every other
result in this paper, a fact about one ray, assembled across
$\Omega$, and nothing here rules it out. \Cref{sec:benedicks}'s
single-ray Benedicks theorem already draws on the same ellipticity of
$L_{-\lambda}$ that underlies it.

It is the other direction, convergence of the integral forcing the
\emph{existence} of a nonzero compactly supported witness $f$, that
the fixed-ray method cannot reach. The witness is built in
\cite{BagchiGangulySarkarThangavelu2022} as a single function on all
of $\mathbb{H}^n$, compactly supported jointly in $z$ and $t$, with a
cutoff depending on $k$ and $\lambda$ together through
$(2k+n)|\lambda|$. By the Paley-Wiener theorem, compact support in
$t$ forces $f^\lambda$, as a function of $\lambda$, to extend to an
entire function of $\lambda$. Being entire, it is then mutually
constrained across every ray at once, so no per-ray construction
glued together over $\lambda$ can produce it. This is the genuine
obstruction. It is not an absence of fixed-ray content in the
hypothesis, but a construction that is irreducibly joint in $(z,t)$,
for a reason with no analogue in any of the existence statements
proved by fixed-ray assembly elsewhere in this paper.
\end{rmk}

\subsection{The joint-ray Benedicks theorem}

The single-ray Benedicks theorem of Section~\ref{sec:benedicks} (Theorem~\ref{thm:benedicks}) fixes $\lambda$ and lets only $k$ vary. Its hypothesis is that, for every $k$, the zero set of $W_{S,\lambda,k}(g)$ has finite measure, and its proof uses only the ellipticity of the single operator $L_{-\lambda}$ at that fixed $\lambda$. Dabra, Dasgupta, and Gulia \cite{DabraDasguptaGulia2026} have recently established a genuinely broader Benedicks-type theorem directly for $\mathcal{F}_S$, with a hypothesis imposed jointly across both $\lambda$ and $k$ simultaneously. This forces $f \equiv 0$ under a single finite-measure condition on a subset $W \subseteq \Omega$ bounded away from the degenerate ray $R_\infty$ (that is, $W \subseteq \{|\lambda|>M\}$ for some $M>0$), rather than on each ray separately.

This joint theorem lies beyond the reach of the method of this paper, though for a different reason than Ingham's theorem above. The hypothesis of \cite{DabraDasguptaGulia2026} is stated directly on subsets $V \subseteq \mathbb{H}^n$ and $W \subseteq \Omega$ of finite measure, not as a condition on $f^{-\lambda}$ ray by ray. Here $W$ is additionally required to stay away from $R_\infty$, since their Hilbert-Schmidt bound on $P_W \circ P_V$ degenerates as $|\lambda|\to 0$. Their proof does not go by fixing $\lambda$ and invoking ellipticity of $L_{-\lambda}$ as Theorem~\ref{thm:benedicks} does. Instead it combines a Hilbert-Schmidt norm bound on the composite projection $P_W\circ P_V$, computed by integrating directly over $\Omega\times\mathbb{C}^n$ and $\mathbb{H}^n$ jointly, with a dilation-based linear-independence argument for the restriction $f(\cdot,0)$ on the full group. Neither step is a single-ray statement that could, in principle, be proved at one fixed $\lambda$ and then assembled. The method is quite different from the ellipticity argument of Section~\ref{sec:benedicks}, not merely a version of it requiring a joint hypothesis. Reaching \cite{DabraDasguptaGulia2026}'s theorem would accordingly mean reproducing their projection-operator argument on $\mathbb{H}^n$ and $\Omega$ directly, which is a different undertaking from anything the fixed-ray dictionary of this paper supplies.

\section{Conclusion}
We looked for an operator through which the Strichartz Fourier
transform $\mathcal{F}_S$ on the Heisenberg group factors on each ray
of the Heisenberg fan $\Omega$. We isolated and named this operator
the Strichartz Weyl transform $W_{S,\lambda}$. We showed it to be, up
to a scalar, the $k$-th spectral projection of the special Hermite
operator. We developed its fixed-ray structural theory in
\Cref{sec:swt}. We put this theory to two uses. In
\Cref{sec:new-results}, we obtained results for which $\mathcal{F}_S$
had no earlier analogue. In \Cref{sec:core-toolkit}, we recovered the
standard foundational toolkit for $\mathcal{F}_S$. This toolkit was
already established in \cite{Thangavelu2023}. We reached it here
almost entirely from idempotency, self-adjointness, and the
Plancherel bound of a single operator. \Cref{sec:hecke-bochner-schwartz}
treated the Hecke-Bochner formula and the image of
$\mathcal{S}(\mathbb{H}^n)$. The latter was reached only by citing the
global theorem of Astengo, Di Blasio, and Ricci. \Cref{sec:limits}
gave a full account of the two results that lie beyond the fixed-ray
method altogether. We gave the precise structural reasons why neither
can be reached ray by ray. The value of the Strichartz Weyl
transform, we suggest, lies less in any single result recovered or
obtained through it. It lies instead in the breadth of what becomes
reachable once $\mathcal{F}_S$'s defining ingredients are recognised
as a decades-old object in disguise. This mirrors how Thangavelu
identified $\mathcal{F}_S$'s own ingredients in the literature that
preceded it.

\end{document}